\documentclass[12pt,reqno]{article}
\def\hybrid{\topmargin 0pt      \oddsidemargin 0pt
        \headheight 0pt \headsep 0pt
        \textwidth 160true mm       
        \textheight 231true mm         
        \marginparwidth 0.0in
        \parskip 3pt plus 1pt   \jot = 1.5ex}
\usepackage{amssymb}
\usepackage{amsmath}
\usepackage{amsthm}

\hybrid
\usepackage{amssymb}
\usepackage{amsmath}
\usepackage{amsthm}
\usepackage{amscd}

\newcommand{\be}[1]{\begin{eqnarray#1}}
\newcommand{\ee}[1]{\end{eqnarray#1}}

\newtheorem{thm}{Theorem}[section]
\newtheorem{lemma}[thm]{Lemma}
\newtheorem{propn}[thm]{Proposition}
\newtheorem{cor}[thm]{Corollary}

\theoremstyle{definition}
\newtheorem{definition}[thm]{Definition}
\newtheorem{example}[thm]{Example}

\theoremstyle{definition}
\newtheorem{remark}[thm]{Remark}

\numberwithin{equation}{section}

\renewcommand{\frak}{\mathfrak}
\renewcommand{\tilde}{\widetilde}

\newcommand{\ot}{\otimes}

\newcommand{\bA}{\mathbb{A}}
\newcommand{\bO}{\mathbb{O}}

\newcommand{\C}{\mathbb{C}}

\newcommand{\R}{\mathbb{R}}
\newcommand{\Ct}{\mathbb{C}[[t]]}
\newcommand{\Cc}{{\mathcal C}}
\newcommand{\Cm}{{\Cc_M^\infty}}
\newcommand{\CM}{{\Cc_M^\infty[[t]]}}
\newcommand{\TT}{\mathcal{T}}
\newcommand{\BB}{\mathcal{B}}
\newcommand{\La}{\Lambda}

\newcommand{\ff}{\varphi}
\newcommand{\WW}{\mathcal{W}}
\newcommand{\Wb}{\mathbf{W}}
\newcommand{\EE}{\mathcal{E}}
\newcommand{\PP}{\mathcal{P}}
\newcommand{\LL}{\mathcal{L}}
\newcommand{\QQ}{\mathcal{Q}}
\newcommand{\FF}{\mathcal{F}}
\newcommand{\Sbf}{\mathbf{S}}
\newcommand{\wi}{\mathbf{w}}
\newcommand{\ad}{\mathrm{ad}}
\newcommand{\ft}{\frac{1}{t}}
\newcommand{\ftt}{\frac{t}{2}}
\newcommand{\OO}{\mathcal{O}}
\newcommand{\cO}{\mathcal{O}}

\newcommand{\g}{\mathfrak{g}}

\newcommand{\tm}{{\frac{1}{t}}}

\newcommand{\Der}{\operatorname{Der}}

\newcommand{\Tm}{\TT^\C_M}
\newcommand{\TM}{\TT^\C_M[[t]]}
\newcommand{\pa}{\partial}
\newcommand{\paa}[1]{\partial/\partial {#1}}
\newcommand{\XX}{{\mathcal X}}
\newcommand{\YY}{{\mathcal Y}}
\newcommand{\De}{\Delta}
\newcommand{\OP}{{\OO_\PP}}

\newcommand{\half}{{\frac{1}{2}}}
\newcommand{\D}{{\mathcal D}}
\newcommand{\cle}{cl_{Ext}}
\newcommand{\clP}{cl_{PSP}}
\newcommand{\clPQ}{{\Bar \tau}}

\newcommand{\Hwt}{\Gamma(M,d\PP_t^\perp)/d(\Gamma(M,\PP_t^\perp))}

\begin{document}
\baselineskip 18pt
\title{Classification of polarized deformation quantizations}

\author
{Joseph Donin\thanks{Partially supported
by Israel Academy of Sciences Grant no. 8007/99-01 }\\
{\normalsize Dept. of Math.,  Bar-Ilan University}}

\date{}
\maketitle

\begin{abstract}
We give a classification of
polarized deformation quantizations on a symplectic manifold
with a (complex) polarization.

Also, we establish a formula which relates the characteristic class
of a polarized deformation quantization to its Fedosov class and
the Chern class of the polarization.
\end{abstract}

\tableofcontents

\section{Introduction}

Let $(M,\omega)$ be a symplectic manifold, $\TT_M$ its
{\em complexified} tangent bundle.
It is known that classes of deformation quantizations on $(M,\omega)$
are in one-to-one correspondence with their Fedosov classes,
the elements
of $\omega+tH^2(M,\Ct)$. The set $\omega+tH^2(M,\Ct)$ may
be interpreted in the following way. Let $\XX$ be the set of formal
closed 2-forms on $M$ of the form $\omega+t\omega_1+t^2\omega_2+\cdots$.
Let $Aut(M)$ be the group of formal automorphisms of $M$ of
the form $e^{tX}$, where $X=X_0+tX_1+t^2X_2+\cdots$, $X_i\in \TT_M$
is a formal vector field.
Since $\omega$ is nondegenerate,
it is easy to see that the orbit of an element $\omega_t\in\XX$
under the action of $Aut(M)$ is $\omega_t+t\cdot d(\Gamma(M,\TT_M^*))$.
It follows from this that $\omega+tH^2(M,\Ct)$ may be identified
with the set of orbits in $\XX$ under the $Aut(M)$-action.
So, the equivalence classes of deformation quantizations on
$(M,\omega)$ are in
one-to-one correspondence with the orbits in $\XX$.

In the paper, we extend that picture
to polarized deformation quantizations (PDQ).

Let $(M,\omega,\PP)$ be a polarized symplectic manifold, i.e.
$\PP$ is a Lagrangian integrable subbundle of the complexified
tangent bundle to $M$. A PDQ on $(M,\omega,\PP)$ is a pair $(\bA_t,\bO_t)$,
where $\bA_t$ is a deformation quantization on $(M,\omega)$ and
$\bO_t$ is a commutative $t$-adically complete subalgebra of $\bA_t$ such
that $\bO_0=\OO_\PP$, the algebra of functions constant along $\PP$.

Let $\YY$ denote the set of pairs $(\omega_t,\PP_t)$, where
$\omega_t\in \XX$ and $\PP_t$ is a polarization of $\omega_t$
such that $\PP_0=\PP$.
Our result is that the equivalence classes of PDQ's on $(M,\omega,\PP)$ are
in one-to-one correspondence
with the orbits in $\YY$ under the $Aut(M)$-action.
Let us describe this correspondence more precisely.

First, we show that any PDQ is equivalent to a polarized
star-product (PSP). By a PSP we mean a triple,
$(\CM,\mu_t,\OO_t)$, where $(\CM,\mu_t)$ is a star-product,
$\OO_t=\OO_{\PP_t}$, the algebra of functions from $\CM$ constant along
a deformed polarization $\PP_t$, and
the multiplication $\mu_t$ satisfies the condition:
$\mu_t(f,g)=fg$ (the usual multiplication) for $f\in\OO_t$, $g\in\CM$.

Further, we assign to any PSP $(\CM,\mu_t,\OO_t)$ a pair
$(\omega_t,\PP_t)\in\YY$ in the following way.
We put $\PP_t=\PP_{\OO_t}$, the sheaf of formal vector fields
annihilating $\OO_t$. The form $\omega_t$ is equal, locally,
to $\sum_idy_i\wedge dx_i$, where $x_i\in\OO_t$, $y_i\in\CM$,
$i=1,...,\half\dim M$, are Darboux coordinates with respect to
$[,]=[,]_{\mu_t}$, the commutator of $\mu_t$; namely
$[x_i,x_j]=[y_i,y_j]=0$, $[y_i,x_j]=\delta_{ij}$.
It turns out that $\omega_t$ is well defined, i.e.
does not depend on the choice of local Darboux coordinates.
Denote the constructed map from PSP's to $\YY$ by $\tau$.

The map $\tau$ turns out to descend to an isomorphism between
the set of classes of PDQ's
and the set $\bigl[\YY\bigr]$ of orbits in $\YY$.
So, we obtain the following commutative diagram of maps:
\be{}\label{funCD}
\begin{CD}
\{\mbox{PSP's}\}&@>\tau>>&\YY\\
@VVV & & @VVV\\
\{\mbox{classes of PDQ's}\}&@>>>&\bigl[\YY\bigr],
\end{CD}
\ee{}
where the left downward arrow is an epimorphism and the bottom arrow
is an isomorphism of sets.

We show that the top arrow $\tau$ is an epimorphism with the
following properties.

1) Two PSP's are equivalent if and only if their images with
respect to $\tau$ lie on the same orbit.

2) Two PSP's $(\CM,\mu_t,\OO_t)$ and
$(\CM,\tilde\mu_t,\tilde\OO_t)$ have the same image with respect
to $\tau$ if and only if $\OO_t=\tilde\OO_t$ and
$[,]_{\mu_t}=[,]_{\tilde\mu_t}$.

3) Let $(\omega_t,\PP_t)=\tau(\CM,\mu_t,\OO_t)$.
Then, the 2-form
\be{}\label{twoform}
\theta_t=\omega_t+\ftt tr(\nabla^2|_{\PP_t})
\ee{}
represents the Fedosov class of the star-product $(\CM,\mu_t)$.
Here $\nabla$ is a connection on $M$ preserving $\omega_t$, $\PP_t$,
and flat on $\PP_t$ along $\PP_t$. We prove that such a connection
always exists and for it $tr(\nabla^2|_{\PP_t})$ belongs to
$\Gamma(M,d\PP^\perp_t)$.

\medskip

By definition,
$\omega_t$ belongs to $\omega+t\Gamma(M,d\PP_t^\perp)$.
So, it follows from (\ref{twoform}) that
the Fedosov class of the star-product $(\CM,\mu_t)$
can be represented by a 2-form belonging to
$\omega+t\Gamma(M,d\PP_t^\perp)$, as well.
In particular, both $(M,\omega_t,\PP_t)$ and
$(M,\theta_t,\PP_t)$ are formal polarized symplectic manifolds
that are deformations of $(M,\omega,\PP)$.

Another consequence of (\ref{twoform}) is the following one.
Let $\bA_t$ be a
deformation quantization on $(M,\omega)$. Suppose its
Fedosov class $cl_F(\bA_t)$ is represented by the
2-form $\theta_t$ that has a polarization $\PP_t$.
Then $\bA_t$ can be extended to a PDQ $(\bA_t,\bO_t)$, i.e.
there exists a commutative subalgebra
$\bO_t\subset \bA_t$ isomorphic to $\OO_{\PP_t}$.

There is the following interpretation of the image $[\omega_t]$ of
element $\omega_t$ in the $\Ct$-module $\Hwt$. Let
$F(\mu_t,\OO_t)=\{a\in\CM; [a,\OO_t]_{\mu_t}\subset\OO_t\}$. Then,
there is the following exact sequence of $\OO_t$-module and Lie
algebra sheaves:
\be{}\label{seqOL}
\begin{CD}
0@>>>\OO_t@>>>F(\mu_t,\OO_t)@>>>Der(\OO_t)@>>>0.
\end{CD}
\ee{}
According to \cite{BB} and \cite{BK}, $F(\mu_t,\OO_t)$ is called
an $\OO_t$-extension of $Der(\OO_t)$. Equivalence classes of such
extensions are described by their extension classes that are
elements of $\Hwt$. We show that $[\omega_t]$ is just the
extension class of (\ref{seqOL}).

Analogously, $-[tr(\nabla^2|_{\PP_t})]$ is the extension class of
the extension
\be{}\label{seqOB}
\begin{CD}
0@>>>\OO_t@>>>\widetilde T_{\det(\PP_t)}@>>>Der(\OO_t)@>>>0,
\end{CD}
\ee{}
where $\widetilde T_{\det(\PP_t})$ is the sheaf of
$\OO_t$-differential operators of order at most one on the
$\OO_t$-line bundle $\det(\PP_t)$.

Note that $-tr(\nabla^2|_{\PP_t})$ divided by $2\pi\sqrt{-1}$
represents the first Chern class of $\PP$, \cite{KN}.
So, formula (\ref{twoform}) gives a relation between the Fedosov
and extension classes of a PDQ.

Among results related to ours we mention the following.

In \cite{RY}, N. Reshetikhin and M. Yakimov
considered the case of a real polarization on $M$
defined by a Lagrangian fiber bundle $M\to B$.

In \cite{Ka1}, A. Karabegov constructed star-products with separation of
variables on K\" ahler manifolds. This case corresponds to two
polarizations on $M$ defined by holomorphic and anti-holomorphic
vector fields.
In the case of quantization on a K\" ahler
manifold with separation of
variables the class of $\omega_t$ in $H^2(M,\Ct)$
coincides with the class defined by Karabegov
in \cite{Ka1}.
A formula relating the Karabegov and Fedosov
classes in case of K\" ahler manifolds
is found in \cite{Ka2}, see also \cite{KS}, \cite{Ne}.

Our proof of the existence of a polarized star-product associated with
any orbit in $\YY$
uses the Fedosov method adapted for the case with
polarization. The analogous method
was applied by M. Bordemann and S. Waldmann, \cite{BW}, for constructing
a quantization with separation of variables on a K\" ahler manifold.

Another approach to proving a formula relating the Fedosov and
extension classes, using the Deligne classes, was presented in
\cite{BD}. Unfortunately, there is a deficiency in the proof of
Lemma 4.3 of that paper relating the extension and Deligne
classes, however, the proof becomes correct for PSP's with the
same polarization.

The paper is organized as follows.

In Section 2, we study cohomologies of the differential
Hochschild complex on $M$ in presence of a distribution.
Also, we prove a version
of the Kostant-Hochschild-Rosenberg
theorem for functions constant along a distribution.
We use these results latter in proving that any
PDQ is equivalent to a PSP.

In Section 3, we introduce a notion
of $\C$-{\em symplectic manifold}, which will be
convenient for our consideration.
This notion is a
generalization of the notion of symplectic manifold.
Namely, we suppose that symplectic form $\omega$
on a $\C$-symplectic manifold
is a {\em complex} one
and, locally, there exist {\em complex} Darboux coordinates
with respect to $\omega$. For a usual symplectic manifold,
when $\omega$ is real, such coordinates
exist by the Darboux theorem.
In this section,
we establish some facts on $\C$-symplectic manifolds
with polarization.
By a polarization of $\omega$ we mean a Lagrangian subbundle, $\PP$,
of the {\em complexified} tangent bundle on $M$ such that, locally on $M$,
there exist Darboux coordinates $x_i, y_i$, $i=1,...,\half\dim M$,
where $x_i\in\OO_\PP$ for all $i$.
So, (pseudo-)K\" ahler manifolds as well as purely real polarizations
are included in our considerations.
Note that from an analog of ``Dolbeault Lemma'' proved in
\cite{Ra} one can derive sufficient conditions for $\PP$ to be
a complex polarization of $\omega$.

In Section 4, we study properties of formal (or deformed)
polarized symplectic manifolds. In Section 5, we prove the
existence of a polarized symplectic
connection on a formal polarized symplectic manifold,
$(M,\omega_t,\PP_t)$,
and with the help of it
introduce the characteristic class $\tilde c_1(M,\omega_t,\PP_t)$
of a polarized symplectic manifold.

In Section 6, we prove some technical statements related to
deformations of Poisson brackets on $M$. Such deformations
appear, in particular, as commutators of star-products.

In Section 7, we study properties of
PDQ's. In particular, we prove the important fact
that any PDQ is equivalent to a PSP.

In Section 8, we define the extension class of a PDQ.
Besides, we assign to any PSP
an element of $\YY$, and to any class of PDQ's
an orbit in $\YY$.
We prove that the later assignment is a monomorphism that, actually,
is an isomorphism, as we show in the next section.

In Section 9, we prove that each element of $\YY$ corresponds to a PSP.
To this end, we adapt the Fedosov method for constructing
a PSP corresponding to a given pair $(\omega_t,\PP_t)\in\YY$.
By this method, a polarized symplectic connection, $\nabla$,
extends to a Fedosov connection on the bundle of Weyl algebras on $M$.
This connection has two scalar curvatures:
the Weyl curvature, $\theta_t$, and the Wick curvature
that turns out to be just $\omega_t$. We show that
these curvatures differ from each other by
$\frac{t}{2}tr(\nabla^2|_{\PP_t})$,
which immediately proves (\ref{twoform}).

In Section 10, we formulate the main theorem
collecting the results of the paper and give some
corollaries.

{\bf Acknowledgments.} I thank J.Bernstein, P.Bressler,
B.Fedosov, A.Karabegov,
and A.Mudrov for helpful discussions.

\section{Complex distributions}
\label{subsec2.2}

For a smooth manifold $M$ we will denote by
$\Cm$ the sheaf of {\em complex valued} smooth functions
on $M$ and by $\Tm=\TT_M\ot_\R\C$ the complexified tangent bundle on $M$.

We say that a set of smooth functions $x_i$, $i=1,...,\dim M$,
given on an open subset $U\subset M$
form a system of (complex) coordinates on $U$,
if $dx_i$ are linearly independent at each point of $U$.

Since any 1-form on $U$ can be uniquely written as $\sum_ia_idx_i$, one can
define vector fields $\pa/\pa x_i\in\Tm$ in the following way.
If $f$ is a function on $U$ and $df=\sum_ia_idx_i$, then
$(\pa/\pa x_i)f=a_i$.

Let $Q$ be a subbundle in a complex vector bundle $E$ over $M$.
We denote by
$Q^\perp$ the subbundle in $E^*$, the complex dual to $E$,
orthogonal to $Q$. If sections $e_i$ form a local frame in $Q$, we
set $(e_i)^\perp=Q^\perp$.

A {\em (complex) distribution} on a manifold $M$ is a subbundle of
$\Tm$.

\begin{definition}
A distribution $\PP$ is said to be {\em integrable} if, locally on $M$,
there exist
(complex valued) functions $f_1,\dots,f_k$ such that $df_1,\dots, df_k$
give a local frame in $\PP^\perp$, i.e. $df_i$ are linearly
independent at each point and $\PP=(df_i)^\perp$.
\end{definition}

An integrable distribution $\PP$ is obviously involutive, i.e.
$[\PP,\PP]\subset\PP$.

Let $\PP$ be an integrable distribution on $M$. We will denote by
$\OP$ the sheaf of functions on $M$ constant along $\PP$, i.e.
$f\in\OP$ if and only if $Xf=0$ for any vector field $X\in\PP$.

\subsection{The Kostant-Hochschild-Rosenberg theorem in
presence of a distribution}

Let $M$ be a smooth manifold.
Let $\D^n$ be the sheaf of $n$-differential operators on $M$
and $\D^\bullet$ the corresponding Hochschild complex with
differential $d$. Let $\wedge^\bullet\TT$ be the complex of sheaves
of polyvector fields on $M$ with zero differential, $\TT=\Tm$.

There is the following ``smooth" version of the
Kostant-Hochschild-Rosenberg theorem, \cite{Ko}, Thm. 4.6.1.1.
\begin{propn}\label{propKHR} The natural embedding
\be{}\label{KHR}
\wedge^\bullet\TT\to\D^\bullet
\ee{}
is a quasiisomorphism of complexes.
Moreover, if $\varphi\in\D^n$ is a Hochschild cocycle, then
its alternation $Alt(\varphi)$ is a polyvector field of $\wedge^\bullet\TT$
cohomological to $\varphi$.
\end{propn}

\begin{proof} Arguments of this proof will be used also in
proving the next proposition.
The proposition is local on $M$, so it is enough
to prove it replacing $M$ by an open set
$U\subset M$ having complex coordinates $x_i$, $i=1,...,\dim M$.
Any differential operator on $U$ may be uniquely presented as
a polynomial in $\paa{x_i}$ with coefficients being smooth
functions on $U$. Hence, $\D^\bullet$ coincides over $U$ with
the complex $\Cc^\bullet(\TT)$.
Here, for any vector bundle $E$, we denote by $\Cc^\bullet(E)$
the complex $(\ot^\bullet Sym(E),d)$ with
differential of the form
$$
d:\ot^n Sym(E)\to \ot^{n+1} Sym(E),
$$
$$
d(a_1\ot...\ot a_n)=1\ot a_1\ot...\ot a_n+
$$
$$
+\sum_{i=1}^n(-1)^ia_1\ot...\ot\De a_i\ot...\ot a_n
+(-1)^{n+1}a_1\ot...\ot a_n,
$$
where $\De$ is the comultiplication in the symmetric algebra
$Sym(E)$ generated by the rule
$\De(a)=a\ot 1+1\ot a$ for $a\in E$.

One has the following well known statement (see, for example,
the proof of Thm. 4.6.1.1. in \cite{Ko}).

\begin{lemma}\label{lemKHR}
Let $E$ be a line bundle over $M$. Then the conclusion
of Proposition \ref{propKHR} holds for the map
\be{}\label{GKHR}
\wedge^\bullet E\to\Cc^\bullet(E).
\ee{}
\end{lemma}

Applying this lemma to $E=\TT$ we prove
the proposition.

\end{proof}

Let $(M,\PP)$ be a smooth manifold with integrable
distribution.
We call an $n$-chain $\nu\in \D^n$ {\em polarized} if
$\nu(a_1,...,a_n)=0$ whenever $a_1,...,a_n\in\OP$.
We call $\nu$ {\em strongly polarized}, if $\nu(a_1,...,a_n)=0$
whenever $a_1,...,a_{n-1}\in\OP$ and $a_n\in\Cm$.

\begin{propn}\label{propdop}
Let $\nu\in \D^2$ be a polarized Hochschild
$2$-cochain such that
$d\nu$ is strongly polarized. Then, there exists a
polarized differential operator $b$ such that
$\nu+db$ is strongly polarized.
\end{propn}

\begin{proof} Since subsheaves of polarized and strongly polarized
cochains are subsheaves of $\Cm$-modules,
it is enough to prove the proposition
locally on $M$. So, let $U$ be an open set with coordinates
$x_i$, $i=1,...,\dim M$ such that $\paa{x_i}$, $i=1,...,k$,
form a local frame in $\PP$.
Let $\QQ$ be the subbundle in $\TT$ generated by $\paa{x_i}$,
$i=k+1,...,n$. Thus, $\TT=\QQ\oplus\PP$ over $U$.
There is the natural isomorphism of complex $\Cc^\bullet(\TT)$
with the tensor product of complexes
$\Cc^\bullet(\QQ,\PP)=\Cc^\bullet(\QQ)\ot\Cc^\bullet(\PP)$,
so we can identify $\Cc^\bullet(\TT)$ with $\Cc^\bullet(\QQ,\PP)$.
Similarly, we identify the complex $\wedge^\bullet(\TT)$
with
$\wedge^\bullet(\QQ,\PP)=\wedge^\bullet(\QQ)\ot\wedge^\bullet(\PP)$.

Thus, the map (\ref{GKHR}) generates the map of
complexes
\be{}\label{gmap}
\wedge^\bullet(\QQ,\PP)\to\Cc^\bullet(\QQ,\PP).
\ee{}
Complex $\Cc^\bullet(\QQ,\PP)$
decomposes obviously into
a direct sum of subcomplexes $\Cc^\bullet_{k,l}(\QQ,\PP)$, where
$\Cc^\bullet_{k,l}(\QQ,\PP)$ consists of elements of total degree
$k$ with respect to $\PP$ and $l$ with respect to $\QQ$.
The same is true for $\La^\bullet(\QQ,\PP)$.
Due to Proposition \ref{propKHR}, the map (\ref{gmap}) is a quasiisomorphism of
bigraded complexes.

It is clear that an element of $\Cc^n(\QQ,\PP)$ is polarized,
if it is a sum of tensor monomials having degree $>0$ with respect to $\PP$.
An element of $\Cc^n(\QQ,\PP)$ is strongly polarized
if it is a sum of tensor monomials
of the form $a_1\ot\cdots\ot a_n$ where $a_1\ot\cdots\ot a_{n-1}$
is polarized.

So, we may suppose that given $\nu$ is a polarized
Hochschild cochain on $U$ belonging to $\Cc^2(\QQ,\PP)$.
It can be written as
$\nu=\nu_0+\nu_1$,
where $\nu_0=\sum (a_i\ot c_i)$, the sum of all
tensor monomials in $\nu$ such that $a_i\in\Cc^1(\QQ)$, $c_i\in\Cc^1(\PP)$.
Let us denote $b=\sum a_ic_i$. It is clear that
$\nu'=\nu+db$ does not contain tensor monomials of that type.
The proposition will be proved if we show that $\nu'$ is
strongly polarized.
Let us prove this.

Suppose $\nu'=\nu'_0+\nu'_1$,
where $\nu'_0$ is not strongly polarized
and $\nu'_1$ is strongly polarized.
Then $\nu'_0$ has the form
\be{}\label{nu}
\nu'_0=\sum (a_i\ot b_i)(1\ot c_i),
\ee{}
where $(a_i\ot b_i)\in \ot^2 Sym(\QQ)$ and
$c_i\in Sym(\PP)$ are linearly independent.
Besides, $b_i$ are of degree $>0$.

The element $d\nu=d\nu'=d\nu'_0+d\nu'_1$ is strongly polarized
by hypothesis of the proposition, $d\nu'_1$ being the coboundary
of a strongly polarized element $\nu'_1$ is strongly polarized too.
All summands of $d\nu'_0$
with first two factors being of degree zero with respect to $\PP$ are
$(d(a_i\ot b_i)+a_i\ot b_i\ot 1)(1\ot 1\ot c_i)$.
These summands are not strongly polarized.
Hence,
$\sum_i (d(a_i\ot b_i)+a_i\ot b_i\ot 1)(1\ot 1\ot c_i)=0$.
Since elements $1\ot 1\ot c_i$ are linearly independent,
it follows that $d(a_i\ot b_i)=-a_i\ot b_i\ot 1$ for all $i$, which
is only possible if $a_i\ot b_i=0$ for all $i$.
It follows from (\ref{nu}) that $\nu'_0=0$. So, $\nu'$ is
equal to $\nu'_1$ which is strongly polarized.
\end{proof}

\subsection{Differential operators in presence of a distribution}

Let $(M,\PP)$ be a smooth manifold with integrable
distribution. Let $z_i$, $y_j$ be
complex coordinates on an open set $U\subset M$
such that $\PP=(dz_i)^\perp$.
Vector fields $\pa/\pa y_j$ form a local frame in $\PP$, since,
by definition, $\OP$ consists of functions $a\in \Cm$ such that
$da$ has the form $\sum_ia_idz_i$. Since $da$ is closed,
$\pa a_i/\pa y_j=0$ for all $i,j$, which implies that
$a_i=\pa a/\pa z_i\in\OP$ for all $i$.

Let $\QQ$ be a subbundle in $\TT=\Tm$ generated by $\paa{z_j}$,
so that $\TT=\QQ\oplus\PP$ over $U$.

The vector bundle $\TT/\PP$ may be considered as the sheaf
of derivations from $\OP$ to $\Cm$, $Der(\OP,\Cm)$.
Locally, such derivations can be presented in the form
$\sum b_i\paa{z_i}$, $b_i\in\Cm$, i.e. as sections of $\QQ$.
Denote by $Der(\OP)$ the $\OP$-submodule of $Der(\OP,\Cm)$
consisting of operators which take $\OP$ to itself.
It is clear that $Der(\OP,\Cm)=\Cm\ot_\OP Der(\OP)$.
Locally, elements of $Der(\OO_\PP)$ have the form
$\sum_ia_i\pa/\pa z_i$, $a_i\in\OP$.

Let $\wedge^\bullet(\TT/\PP)$
denote the complex of sheaves of polyvector fields on $M$ from
$\OP$ to $\Cm$.

Let $\D^\bullet(\OP,\Cm)$ denote
the restriction of the Hochschild complex $\D^\bullet$ to $\OP$.
So, the sheaf $\D^n(\OP,\Cm)$ may be considered as the sheaf
of $n$-differential operators from $\OP$ to $\Cm$. Locally,
elements of $\D^1(\OP,\Cm)$, the sheaf of differential
operators from $\OP$
to $\Cm$, can be presented as polynomials
in $\paa{z_i}$ with smooth coefficients. So, locally on $M$,
complex $\D^\bullet(\OP,\Cm)$ is isomorphic to the complex
$\Cc^\bullet(\QQ)$ (see previous subsection).

We will need the following version of the Kostant-Hochschild-Rosenberg
theorem.

\begin{propn}\label{propKHR1}
The natural embedding
\be{}\label{imbO}
\wedge^\bullet(\TT/\PP)\to\D^\bullet(\OP,\Cm)
\ee{}
is a quasiisomorphism of complexes. Moreover,
if $\ff\in\D^n(\OP,\Cm)$ is a cocycle, then $Alt(\ff)$
is a polyvector field of $\wedge^n(\TT/\PP)$
cohomological to $\ff$.
\end{propn}

\begin{proof}
As follows from above, embedding (\ref{imbO})
is locally isomorphic to the embedding
\be{*}
\wedge^\bullet\QQ\to\Cc^\bullet(\QQ).
\ee{*}
Now the proposition follows from Lemma \ref{lemKHR} when $E=\QQ$.
\end{proof}

\begin{remark}
All conclusions of Propositions \ref{propKHR}, \ref{propdop},
\ref{propKHR1} remain true for global sections of
the corresponding sheaves, since they are
sheaves of $\Cm$-modules.
\end{remark}

\subsection{Differential forms in presence of a distribution}

Let $(M,\PP)$ be a smooth manifold with integrable distribution.

The sheaf $\PP^\perp=(\TT/\PP)^*$ of differential forms on $M$
which being applied to vector fields from $\PP$ give zero,
may be written as $\Cm d\OO_\PP$.

Denote $\Omega_{\OO_\PP}^1=Hom_{\OO_\PP}(Der(\OP),\OP)$,
the sheaf of 1-forms on $\OP$. It is clear that
$\Omega_{\OO_\PP}^1=\OP d\OP$.

Denote by $\Omega_\OP^{1cl}$ the subsheaf of closed forms
of $\Omega_\OP^1$.

\begin{lemma}\label{lemma2.2}

a) The sequence of sheaves
\be{}\label{seqomega}
\begin{CD}
0 @>>> \C @>>> \OP @>d>> \Omega_\OP^{1cl} @>>> 0
\end{CD}
\ee{}
is exact.

b) The sequence of sheaves
\be{}\label{seqphi}
\begin{CD}
0 @>>> \Omega_\OP^{1cl} @>>> \PP^\perp @>d>> d\PP^\perp @>>> 0
\end{CD}
\ee{}
is exact.
\end{lemma}

\begin{proof}
Let functions $z_i$, $y_j$ form a local basis on $M$ and
$\PP=(dz_i)^\perp$. Let us prove a). It is sufficient to establish
the exactness at the third term of (\ref{seqomega}).
Let $\alpha=\sum_i a_idz_i\in\Omega_\OP^{1cl}$. Since
$\alpha$ is closed, there exists, locally, $f\in\Cm$ such that
$df=\alpha$. Since $\alpha$ does not contain terms of
the form $gdy_j$, one has $\pa f/\pa y_j=0$ for all $j$.
Hence, $f\in\OO_\PP$.

To prove b), it is sufficient to establish
the exactness at $\PP^\perp$. To this end, suppose
$\beta=\sum_i b_idz_i\in\PP^\perp$ and $d\beta=0$.
Since $\beta$ is closed,
$\pa b_i/\pa y_j=0$ for all $i,j$. This means that all $b_i\in\OO_\PP$.
Thus, $\beta\in\Omega_\OP^1$ and closed, i.e. $\beta\in\Omega_\OP^{1cl}$.
\end{proof}

\begin{propn}\label{fundiso}
Let $(M,\PP)$ be a smooth manifold with an integrable distribution.
Then, there is the natural isomorphism
\be{}\label{fundiso1}
H^1(M,\Omega_\OP^{1cl})\backsimeq \Gamma(M,d\PP^\perp)/d(\Gamma(M,\PP^\perp)).
\ee{}
\end{propn}

\begin{proof} This is an immediate consequence of the cohomological
exact sequence
for (\ref{seqphi}) and of $H^i(M,\PP^\perp)=0$ for $i>0$.
\end{proof}

\section{$\C$-symplectic manifolds and their polarizations}
\subsection{$\C$-symplectic manifolds}

\begin{definition} By a {\em $\C$-symplectic manifold}
we mean a pair $(M,\omega)$,
where $M$ is a smooth manifold
and $\omega$ a closed nondegenerate {\em complex} 2-form on $M$
satisfying the following condition:
each point of $M$ has a neighborhood $U$ and complex coordinates
$x_i,y_i$, $i=1,...,\half\dim M$, on $U$
such that the form $\omega$ on $U$ can be presented as
\be{}
\omega=\sum_idy_i\wedge dx_i,
\ee{}
\end{definition}
If the form $\omega$ is real, such a presentation is possible
by the Darboux theorem.
One has
\be{}\label{darb}
\{x_i,x_j\}=\{y_i,y_j\}=0, \ \ \
\{y_i,x_j\}=\delta_{ij}
\ee{}
for all $i,j$, where $\{\cdot,\cdot\}$ is the Poisson bracket
inverse to $\omega$.
We call such functions $x_i,y_i$ {\em Darboux coordinates}
with respect to $\omega$ (or $\{\cdot,\cdot\}$).

In what follows we only deal with $\C$-symplectic
manifolds, so we simply call them symplectic ones.

\subsection{Polarization}

\begin{definition}\label{def1.1}
Let $(M,\omega)$ be a symplectic manifold.
We call a (complex) distribution $\PP$ on $M$
a polarization of $\omega$, if, locally,
there exist (complex) Darboux coordinates,
$x_i$, $y_i$, with respect to $\omega$
such that $\PP=(dx_i)^\perp$, i.e. $x_i\in\OO_\PP$.

We call the triple $(M,\omega,\PP)$
a {\em polarized symplectic manifold} (PSM).
\end{definition}

It follows that a polarization of $\omega$ is, in particular,
an integrable distribution and a Lagrangian subbundle
with respect to $\omega$.

\begin{propn}
Let $(M,\omega,\PP)$ be a PSM. Then
$\omega\in\Gamma(M,d\PP^\perp)$.
\end{propn}

\begin{proof}
Let $x_i$, $y_i$ be local Darboux coordinates on $M$
such that $\PP=(dx_i)^\perp$ and $\omega=\sum dy_i\wedge dx_i$.
Then, locally, $\omega=d(\sum y_i dx_i)$ and
$\sum y_i dx_i\in\PP^\perp$.
\end{proof}

\begin{propn}\label{propmax} Let $(M,\omega,\PP)$ be
a PSM. Then, $\OO_\PP$ is a maximal
commutative Lie subalgebra in $\Cm$ with respect to
the Poisson bracket $\omega^{-1}$.
\end{propn}

\begin{proof}
It follows from Definition \ref{def1.1} that, locally,
the bracket $\omega^{-1}$ may be written in the form
\be{}\label{eq10}
\{\cdot,\cdot\}=\sum_i\bar\partial_i\wedge\partial_i,
\ee{}
where $\partial_i=\{y_i,\cdot\}$, $\bar\partial_i=\{\cdot,x_i\}$.
The module $\OO_\PP$ consists, locally, of elements $g\in\Cm$ such
that $\{g,x_i\}=\bar\partial_ig=0$ for all $i$.
Putting two such elements $g_1,g_2$ in (\ref{eq10}),
we obtain that $\{g_1,g_2\}=0$. So $\OO_\PP$ is commutative.
The maximality of $\OO_\PP$ is obvious. Indeed, if $a\in\Cm$ commutes
with $\OO_\PP$, then, in particular, $\{x_i,a\}=0$ for all $i$, hence
$a\in\OO_\PP$.
\end{proof}

\section{Deformations of a polarized symplectic manifold}

\subsection{Formal everything}

Let $t$ be a formal parameter. We will consider on $M$ formal
functions, formal vector fields, formal forms, etc.,
which are elements of $\CM$, $\TM$, $\Phi^k[[t]]$, etc.
In the formal case all sheaves over $M$ and their morphisms
will be sheaves and morphisms of $\Ct$-modules.

Let $\BB$ be a sheaf over $M$. We call the map
$\sigma:\BB[[t]]\to \BB$, $b_0+tb_1+\cdots\to b_0$,
the {\em symbol map}. For a subsheaf $\FF_t\subset\BB[[t]]$,
we denote $\FF_0=\sigma(\FF_t)\subset\BB$.

Let $\FF_t$ be a subsheaf of $\BB[[t]]$. We call $\FF_t$
a $t$-{\em regular} subsheaf if it is complete in $t$-adic
topology, and $tb\in\FF_t$, $b\in\BB[[t]]$ implies
$b\in\FF_t$.

For a $t$-regular subsheaf $\FF_t\subset\BB[[t]]$ the natural
map $\FF_t/t\FF_t\to\FF_0$ is an isomorphism. In this case we call $\FF_t$
a {\em deformation} of $\FF_0$.

Let $\BB$ be a vector bundle over $M$. Then, a subsheaf
$\FF_t\subset\BB[[t]]$ is called a (formal) subbundle,
if it is a $t$-regular subsheaf of $\CM$-modules.

All notions and statements above carry over to the formal case.

For example, a deformation of a distribution $\PP$ on $M$ is
a subbundle, $\PP_t$, of $\TM$ such that $\PP_0=\PP$.
A system of (formal) coordinates on an open set $U\subset M$
is a set of formal functions $x_i=x_{i0}+tx_{1i}+\cdots$, $i=1,...,\dim M$,
$\Ct$-linearly independent at each point of $U$.
A formal polarized symplectic manifold is a triple
$(M,\omega_t,\PP_t)$, where $\omega_t$ is a formal symplectic
form on $M$
(i.e. a closed 2-form of the form
$\omega_t=\omega_0+t\omega_1+\cdots$
with nondegenerate $\omega_0$ )
and $\PP_t$ is a polarization of $\omega_t$ in sense of
Definition \ref{def1.1}, i.e.,
locally, there exist formal Darboux coordinates $x_i$, $y_i$ with
respect to $\omega_t$ such that $\PP_t=(dx_i)^\perp$.

We say that $(M,\omega_t,\PP_t)$ is
a deformation of a polarized symplectic manifold $(M,\omega,\PP)$,
if $(M,\omega_0,\PP_0)=(M,\omega,\PP)$.

The following proposition, which follows from
Propositions \ref{prop12} and \ref{prop11} below,
shows that any deformation of (polarized) symplectic
manifold is a formal (polarized) symplectic manifold.

\begin{propn} a) Let $\omega_t$ be a formal closed 2-form.
If $\omega_0$ admits, locally,
Darboux coordinates,
then there exist their lifts being formal Darboux coordinates
for $\omega_t$.

b) Let $\PP_t$ be an integrable distribution
and Lagrangian with respect to $\omega_t$.
Let $\PP_0$ be a polarization of $\omega_0$. Then
$\PP_t$ is a polarization of $\omega_t$.
\end{propn}

\subsection{Local structure of
deformed polarizations}

It is clear that formal vector fields $tX_1+t^2X_2+\cdots\in t\TM$
form a sheaf of pro-nilpotent Lie algebras. It follows
that elements $e^{tX}$, $X\in \TM$, form a sheaf of pro-unipotent
Lie groups of formal automorphisms of $M$.

Let $x_i$ be formal coordinates on $U$ and $a_i$, $i=1,...,\dim M$,
arbitrary formal functions on $U$. Then there exists a derivation, $D$,
of $\CM$ that takes $x_i$ to $a_i$. Such a derivation is
$D=\sum_i a_i\pa/\pa x_i$. This implies the following

\begin{lemma}\label{lemma1.1}
Let $x_i$, $x'_i$, $i=1,...,\dim M$, be two systems of
formal coordinates on an open set
$U\subset M$, and $x_i=x'_i \mod t$. Then, there exists a formal
automorphism on $U$ that takes $x_i$ to $x'_i$.
\end{lemma}

\begin{propn}\label{prop19}
a) Let $\PP_t$ be an integrable distribution on $M$ which is a deformation
of a distribution $\PP$.
Then, locally, there exists a formal vector field, $X$,
such that $e^{tX}$ gives an isomorphism $\PP_t$ with $\PP[[t]]$.

b) Let $(M,\omega_t,\PP_t)$ be a deformation of a polarized symplectic
manifold $(M,\omega,\PP)$. Then,
each point of $M$ has a neighborhood, $U$,
and a formal vector field $X$ on $U$
such that $e^{tX}$ gives an isomorphism of
$(M,\omega_t,\PP_t)|_U$ with
the trivial deformation $(U,\omega,\PP[[t]])$.
\end{propn}

\begin{proof} a) Locally, there exist functions $x_{it}=x_{i0}+tx_{i1}+\cdots$,
$i=1,...,k$, such that $\PP_t=(dx_{it})^\perp$ and hence $\PP=(dx_{i0})^\perp$.
Let us add functions $x_{(k+1)0},...,x_{n0}$ in such a way that
all $x_{i0}$, $i=1,...,n$, form a coordinate system.
According to Lemma \ref{lemma1.1}, there exists a formal automorphism
which takes coordinates $x_{it}$, $i=1,...,k$, $x_{j0}$, $j=k+1,...,n$,
to coordinates $x_{i0}$, $i=1,...,n$. This formal automorphism
gives obviously an isomorphism $\PP_t$ onto $\PP[[t]]$.

b) Let $U$ be an open set in $M$ where Darboux coordinates
$x_{it}=x_{i0}+tx_{i1}+\cdots$, $y_{it}=y_{i0}+tx_{i1}+\cdots$,
$i=1,...,\half\dim M$, exist, and $\PP_t=(dx_{it})^\perp$.
Then, $x_{i0}$, $y_{i0}$ are Darboux coordinates with respect
to $\omega_0=\omega$ such that $\PP=(dx_{i0})^\perp$.
By Lemma \ref{lemma1.1} there exists a formal vector field $X$ on $U$
such that the formal automorphism $e^{tX}$ takes coordinates
$x_{it}$, $y_{it}$ to $x_{i0}$, $y_{i0}$. Such $X$ satisfies the
conclusion of the proposition.
\end{proof}

Let $\PP$ be an integrable distribution. Denote by $\OO_\PP$
the sheaf of functions constant along $\PP$. Let $\PP_t$ be
a deformation of $\PP$. It follows from the previous proposition
that, locally, the pair $(\CM,\OO_{\PP_t})$ is isomorphic to
the pair $(\CM,\OO_\PP[[t]])$,
hence $\OO_{\PP_t}$ is a $t$-regular subalgebra of $\CM$.
One holds the following inverse assertion.

\begin{propn}\label{propdis} Let $\PP$ be an integrable distribution on $M$.
Let $\OO_t$ be a $t$-regular subalgebra of $\CM$ such that
$\OO_0=\OO_\PP$. Then, there exists a deformation, $\PP_t$, of $\PP$
such that $\OO_t=\OO_{\PP_t}$.
\end{propn}

\begin{proof} Let us prove that
$d\OO_t\subset (\Tm)^*[[t]]$ is $t$-regular.
Let $b\in(\Tm)^*[[t]]$
and $tb\in d\OO_t$. Then, there exists $a=a_0+ta'\in\OO_t$ such
that $da=tb$. It follows that $a_0$ is a constant, so $ta'\in\OO_t$.
Since $\OO_t$ is $t$-regular, $a'\in\OO_t$, too.
Therefore, $b=da'\in d\OO_t$, so that $d\OO_t$
is a $t$-regular submodule in $(\Tm)^*[[t]]$.
This implies that $\CM d\OO_t$ is a subbundle in $(\Tm)^*[[t]]$,
so $\PP_t=(d\OO_t)^\perp$ is a subbundle in $\TM$.
Moreover, $\PP_0=\PP$. So, $\PP_t$ is a deformation of $\PP$.
One has $\OO_t\subset\OO_{\PP_t}$. Since these two subalgebras
are $t$-regular and
coincide at $t=0$, we have
$\OO_t=\OO_{\PP_t}$.
\end{proof}

The last proposition shows that there is a one-to-one correspondence
between deformations of an integrable distribution $\PP$
and deformations of $\OO_\PP$.

\subsection{Action of formal automorphisms on a symplectic form}

Let $(M,\omega_t)$ be a formal symplectic manifold, $\TT=\Tm$.

The well known formula for the Lie derivative
$$L_X=i(X)d+di(X)$$
implies
\be{}
L_X\omega_t=d\alpha(X),
\ee{}
where $\alpha: \TT[[t]]\to \TT^*[[t]]$ is the map defined by
$X\mapsto\omega_t(X,\cdot)$.
Since $\alpha$ is an isomorphism,
we have the following lemma.
\begin{lemma}\label{lemma1}
The orbit of $\omega_t$ by the action of the group of formal
automorphisms $e^{tX}$, $X\in\Gamma(M,\TT[[t]])$, is
$\omega_t+td(\Gamma(M,\TT^*[[t]]))$.
\end{lemma}
The lemma shows that the orbit of $\omega_t$ coincides with
the cohomology class of $\omega_t$ in $\omega_0+tH^2(M,\Ct)$.

Let $(M,\omega_t,\PP_t)$ be a formal PSM and
$X\in\Gamma(M,\PP_t)$. Since $\alpha(\PP_t)=\PP_t^\perp$,
one has $L_X\omega_t=d\alpha(X)\in \Gamma(M,d\PP_t^\perp)$.
The argument as above implies

\begin{lemma}\label{lemma3}
By the action of the group generated
by $e^{tX}$, $X\in\Gamma(M,\PP_t)$,
the orbit of $\omega_t$   is
$\omega_t+td(\Gamma(M,\PP_t^\perp))$.
\end{lemma}
The lemma shows that the orbit
of $\omega_t$ corresponds to the cohomology class of
$\omega_t$ in $\Hwt$.

\section{Polarized symplectic connection
and characteristic class of a polarized symplectic manifold}\label{PSC}

\subsection{Polarized symplectic connection}

Let $(M,\omega,\PP)$ be a (formal) PSM. Denote, for shortness, $\TT=\TM$.

\begin{definition}
We call a connection, $\nabla$, on $M$ a $\PP$-{\em symplectic connection} if

a) it preserves $\omega$ and is torsion free, i.e. is a symplectic connection;

b) it preserves $\PP$, i.e. $\nabla_X(\PP)\subset\PP$ for any $X\in\TT$;

c) it is flat on $\PP$ along $\PP$, i.e. for any $X,Y\in\PP$ one has
$(\nabla_X\nabla_Y-\nabla_Y\nabla_X-\nabla_{[X,Y]})(\PP)=0$.
\end{definition}

\begin{propn}\label{prop1.6}
Let $(M,\omega,\PP)$ be a (formal) PSM. Then, there exists a
$\PP$-symplectic connection on $M$.
\end{propn}

\begin{proof}
Let functions $a_1,...,a_{2n}$, $2n=\dim M$,
form local Darboux coordinates on an open set $U\subset M$ and be such that
$a_i\in\OO_\PP$ for $i=1,...,n$.
Let $X_i=X_{a_i}$ be the corresponding Hamiltonian vector fields.
Then, vector fields $X_i$, $i=1,...,n$,
form a local frame in
$\PP$. Also, all $X_i$ commute and form a local frame in $\TT$.
Let $\nabla$ be the standard flat connection on $U$ associated with coordinates $a_i$.
This connection is defined on $U$ by the rule $\nabla_{X_i}X_j=0$.

It is easy to see that $\nabla$ is a $\PP$-symplectic connection on $U$.
Moreover,
since $X_f\in\PP$ is equivalent to $df\in\PP^\perp$, the connection
$\nabla$ satisfies the following property for Hamiltonian
vector fields:
\be{}\label{hamilt}
\nabla_{X_f} X_g=0 \quad\quad \mbox{for\ } X_f,Y_g\in\PP.
\ee{}
Now, let us prove the existence of a global connection.

Let $\{U_\alpha\}$ is an open covering of $M$ such that on each $U_\alpha$
there is a $\PP$-symplectic connection $\nabla_\alpha$ as above.
Then $\nabla_\alpha-\nabla_\beta$ defined on
$U_\alpha\cap U_\beta$ form a \v Cech
cocycle $\psi_{\alpha,\beta}\in Hom(\TT\ot\TT, \TT)$,
$\psi_{\alpha,\beta}(X,Y)=\nabla_{\alpha X}Y-\nabla_{\beta X}Y$.

Elements $\psi_{\alpha,\beta}$ satisfy the following properties.

As follows from (\ref{hamilt}),
\be{}\label{hamilt1}
\psi_{\alpha,\beta}(X,Y)=0 \quad\quad \mbox{for\ }
X,Y\in\PP.
\ee{}
Since all $\nabla_\alpha$ are torsion free, $\psi_{\alpha,\beta}$ are symmetric. Since all
$\nabla_\alpha$ preserve $\PP$, one has $\psi_{\alpha,\beta}(X,Y)\in \PP$ for $Y\in\PP$.

In addition, $\psi_{\alpha,\beta}$ considered as elements
from $Hom(\TT, Hom(\TT,\TT))$,
$X\mapsto\psi_{\alpha,\beta}(X,\cdot)$, belong to
$Hom(\TT,\frak{sp}(\TT))$, where $\frak{sp}(\TT)$ consists of
endomorphisms of $\TT$ preserving $\omega$.

Since all the properties above are $\CM$-linear,
one can find tensors $\psi_\alpha\in
Hom(\TT\ot\TT, \TT)$ satisfying all of them and such that
$\psi_\alpha-\psi_\beta=\psi_{\alpha,\beta}$.
Then $\nabla=\nabla_\alpha-\psi_\alpha=\nabla_\beta-\psi_\beta$
is a globally defined connection.
Flatness of $\nabla$ on $\PP$ along $\PP$ follows from the fact that
for all $\alpha$ tensors $\psi_\alpha$ satisfy property (\ref{hamilt1}), i.e.
$\psi_\alpha(X,Y)=0$ for $X,Y\in\PP$.
Also, $\nabla$ is torsion
free because all $\psi_\alpha$ are symmetric. So, $\nabla$ satisfies the proposition.
\end{proof}

\subsection{Characteristic class of a polarized symplectic manifold}
\label{subsec5.2}

Let $\nabla$ be a $\PP_t$-symplectic connection on
a formal polarized symplectic manifold $(M,\omega_t,\PP_t)$.
Let us denote by $\nabla^2|_{\PP_t}$ the curvature
of $\nabla$ restricted to $\PP_t$. Then, $tr(\nabla^2|_{\PP_t})$
is a closed 2-form on $M$ that represents, up to a constant factor,
the first Chern class of $\PP_0$.

\begin{lemma}\label{lempolform}
Let $\nabla$ be a $\PP_t$-symplectic connection.
Then $tr(\nabla^2|_{\PP_t})$ belongs to
$\Gamma(M,d\PP^\perp_t)$.
If $\nabla_1$ is another
$\PP_t$-symplectic connection, then
$tr(\nabla_1^2|_{\PP_t})$ differs from $tr(\nabla^2|_{\PP_t})$
by an element of $d(\Gamma(M,\PP^\perp_t))$.
\end{lemma}

\begin{proof} Follows from flatness of $\PP_t$
along $\PP_t$ with respect to the connection.
\end{proof}
The lemma allows us to consider
the element of $\Hwt$
represented by the form $tr(\nabla^2|_{\PP_t})$,
where $\nabla$ is a
$\PP_t$-symplectic connection, as a
characteristic class of the polarized symplectic manifold
$(M,\omega_t,\PP_t)$.

Due to Proposition \ref{fundiso}, one can also consider
this class
as an element of $H^1(M,\Omega^{1cl}_{\OP_t})$.

\section{Deformations of Poisson brackets}

In this section we prove three technical statements which
we use through the paper.

Let $\pi_0=\{\cdot,\cdot\}$ be a nondegenerate Poisson bracket
on a smooth manifold $M$ of dimension $2n$.
We say that a formal sum $\pi_t=\pi_0+t\pi_1+\cdots$
is a deformation of
$\pi_0$ if all $\pi_i$ are {\em bidifferential operators} on $M$ and
$\pi_t$ defines a Lie algebra structure on the sheaf $\CM$.
We will also denote $\pi_t$ by $[\cdot,\cdot]$. Let us recall the symbol map
$\sigma:\CM\to\Cm$, $a=a_0+ta_1+\cdots\mapsto a_0$. We call $a$ a lift of $a_0$.
We say that functions $\hat x_i$, $\hat\xi_i$, $i=1,...,n$ on an open set
$U\subset M$ form Darboux coordinates with respect to
$[\cdot,\cdot]$, if
$[\hat x_j,\hat x_k] = [\hat\xi_j,\hat\xi_k] = 0$,
$[\hat\xi_j,\hat x_k] = \delta_{jk}$ for all $j,k$.
It is clear that then functions $x_i=\sigma(\hat x_i)$,
$\xi_i=\sigma(\hat\xi_i)$ form Darboux coordinates
with respect to $\{\cdot,\cdot\}$.

\begin{propn}\label{prop12}
Let $[\cdot,\cdot]$ be a deformation of a Poisson bracket
$\{\cdot,\cdot\}$ on $M$.
Let
$\hat x^{(i)}_1\ldots,\hat x^{(i)}_n,
\hat\xi^{(i)}_1,\ldots,\hat\xi^{(i)}_n\in\CM$,
$i=1,2$,
be two systems of Darboux coordinates with respect to $[\cdot,\cdot]$
on a contractible open set $U$ in $M$ satisfying
\be{*}
\sigma(\hat x^{(1)}_j) = \sigma(\hat x^{(2)}_j), \ \
\sigma(\hat\xi^{(1)}_j) = \sigma(\hat\xi^{(2)}_j)
\ee{*}
Then, there exists $B\in \CM$ on $U$ such that the automorphism
$\Phi=\exp(t\cdot ad(B))$,
where $ad(B)=[B,\cdot]$,
satisfies $\Phi(\hat x^{(1)}_j) = \hat x^{(2)}_j$ and
$\Phi(\hat\xi^{(1)}_j) = \hat\xi^{(2)}_j$.
\end{propn}

\begin{proof}
Let $B_0 = 0$. Assume that $B_m$ is such that the automorphism
$\Phi_m = \exp(t\cdot ad(B_m))$
satisfies the conclusion of the proposition modulo $t^{m+1}$.
This assumption
is valid for $m=0$.

Then,
$\hat x^{(2)}_j = \Phi_m(\hat x^{(1)}_j) + t^{m+1}y_j \mod t^{m+2}$,
$\hat\xi^{(2)}_j =\Phi_m(\hat\xi^{(1)}_j) + t^{m+1}\eta_j \mod t^{m+2}$
for suitable
$y_j, \eta_j\in\Cm$. The Darboux commutation relations for
$\hat x^{(i)}_1\ldots,
\hat x^{(i)}_n,\hat\xi^{(i)}_1,\ldots,\hat\xi^{(i)}_n$
imply that the functions $y_1,\ldots,y_n,\eta_1\ldots,
\eta_n$ satisfy
$\{x_j,y_k\}-\{x_k,y_j\}=0$, $\{\xi_j,\eta_k\}-\{\xi_k,\eta_j\}=0$,
$\{\xi_j,y_k\}-\{x_k,\eta_j\} =0$,
where $x_j=\sigma(\hat x^{(1)}_j) = \sigma(\hat x^{(2)}_j)$,
$\xi_j=\sigma(\hat\xi^{(1)}_j) = \sigma(\hat\xi^{(2)}_j)$.
Equivalently, the differential form $\alpha =\sum_j y_jdx_j + \eta_jd\xi_j$ is
closed. By Poincar\'e Lemma there exists $f\in\Cm$ such that $df=\alpha$,
equivalently $y_j=\{\xi_j,f\}$, and $\eta_j=\{f,x_k\}$.

There exists $B_{m+1}\in\Cm$ such that
\be{*}
\exp(t\cdot ad(B_{m+1}) = \exp(ad(t^{m+1} f))\circ\exp(t\cdot ad(B_m))
\ee{*}
and $B_{m+1} = B_m\mod t^{m+1}$. The limit $B = \underset{m\to\infty}{\lim}B_m$
exists and satisfies the conclusions of the proposition.
\end{proof}

\begin{propn}\label{prop11}
Let $[\cdot,\cdot]$ be a deformation of a Poisson bracket
$\{\cdot,\cdot\}$ on $M$ and $\OO_t$ a $t$-adically
complete submodule in $\CM$
being a commutative
Lie subalgebra with respect to $[\cdot,\cdot]$.
Let functions $x_i\in\OO_0$, $\xi_i\in\Cm$ form Darboux
coordinates with respect to $\{\cdot,\cdot\}$ on a contractible
open set $U\subset M$.
Then, there exist their lifts $\hat x_i\in\OO_t$, $\hat\xi_i\in\Cm$
on $U$ which are
Darboux coordinates with respect to $[\cdot,\cdot]$.
\end{propn}

\begin{proof}
Since, by definition, $\OO_t\to\sigma(\OO_t)$ is surjective,
we choose arbitrary lifts $\hat x_j\in \OO_t$ of $x_j$.
Note that $\hat x_j$, $\xi_j$ satisfy the both conclusions modulo $t$.

Let $\hat x_{j,1}:=\hat x_j$, $\hat \xi_{j,1} :=\xi_j$.
Suppose that $m\geq 2$ and $\hat x_j\in\OO_t$, $\hat\xi_{j,m}\in\CM$
satisfies the conclusion of
the proposition modulo $t^{m+1}$.
The assumption on $x_{j,m}$, $\xi_{j,m}$ implies that
\be{*}
[\hat x_j,\hat x_k]&=&0, \\ \
[\hat\xi_{j,m},\hat x_k]& =& \delta_{jk} + y_{jk}t^{m+1}\mod t^{m+2},\\ \
[\hat \xi_j,\hat \xi_k]&=&t^{m+1}z_{jk}\mod t^{m+2}
\ee{*}
for suitable $y_{jk}, z_{jk}\in\Cm$.
The Jacobi identity and the commutation relations for
$x_j$, $\xi_j$ imply that
\be{*}
\{y_{jk},x_l\}-\{y_{jl},x_k\}=0,\\ \
\{z_{lj},x_k\}+\{\xi_j,y_{lk}\}-\{\xi_l,y_{jk}\}=0,\\ \
\{z_{jk},\xi_l\} + \{z_{kl},\xi_j\} + \{z_{lj},\xi_k\} = 0.
\ee{*}
These identities say that the differential form
$$\alpha=y_{jk}dx_j\wedge d\xi_k+z_{jk}dx_j\wedge dx_k$$
is closed.

By the Poincar\'e Lemma there exists a 1-form
$\beta=\sum a_id\xi_i+\sum b_idx_i$ on $U$ such that $d\beta=\alpha$.

Note that $\{x_k,a_i\}=0$ for all pairs $i,k$, since
$\alpha$ contains no terms of the form $fd\xi_i\wedge d\xi_k$.
Hence, $a_i\in\sigma(\OO_t)$. Let $\hat a_i$ be a lift
of $a_i$ in $\OO_t$.

It is easy to check that $\hat x_{i,m+1}=x_{i,m}+t^{m+1}\hat a_i$
and $\hat\xi_{i,m+1}+t^{m+1}b_i$ satisfy the conclusions of
the proposition modulo $t^{m+2}$. Hence, the limits
$\hat x_i = \underset{m\to\infty}{\lim}\hat x_{i,m}$,
$\hat \xi_i= \underset{m\to\infty}{\lim}\hat \xi_{i,m}$
exist and satisfy the conclusions of the proposition.
\end{proof}

We will need a more strong assertion.

\begin{propn}\label{prop11a}
Let $[\cdot,\cdot]$ be a deformation of a Poisson bracket
$\{\cdot,\cdot\}$ on $M$ and $\OO_t$ a
$t$-regular submodule in $\CM$, which is a commutative
Lie subalgebra with respect to $[\cdot,\cdot]$.
Let functions $x_i\in\OO_t$, $\xi_i\in\CM$ form Darboux
coordinates modulo $t^k$, $k>0$, with respect to
$[\cdot,\cdot]$ on a contractible
open set $U\subset M$.
Then, there exist functions $a_i\in\OO_t$, $b_i\in\CM$ on $U$
such that the functions $x_i+t^ka^i$, $\xi_i+t^kb_i$
form Darboux coordinates with respect to $[\cdot,\cdot]$.
\end{propn}
\begin{proof} The same as of Proposition \ref{prop11}.
\end{proof}

\section{Deformation quantization on a polarized symplectic manifold}

\subsection{Deformation quantization}
Let us recall some definitions and facts about
the deformation quantization on a smooth manifold $M$,
see \cite{BFFLS}, \cite{De}, \cite{Fe}.

\begin{definition} a) Let $\Cm$ be the sheaf of smooth
complex valued functions
on $M$. A {\em formal deformation} of $\Cm$ is a sheaf
of $\C[[t]]$-algebras, $\bA_t$, with an epimorphism
$\sigma:\bA_t\to \Cm$ of $\Ct$-algebras
(called  the {\em symbol map})
satisfying the condition: There exists an isomorphism
of $\Ct$-modules $\bA_t\to\CM$ commuting with symbol maps.
(Recall that the symbol map $\sigma:\CM\to\Cm$ takes
$f_0+tf_1+\cdots$ to $f_0$.)

b) Two formal deformations $\bA_{t1}$ and $\bA_{t2}$ of $\Cm$ are
equivalent if there exists a map of sheaves of $\Ct$-algebras
$\bA_{t1}\to\bA_{t2}$, commuting with symbol maps.
\end{definition}

Suppose $\bA_t$ is a formal deformation of $\Cm$. The formula
\be{}
\{a,b\}=\sigma(\frac{1}{t}[\tilde a,\tilde b]),
\ee{}
where $a$ and $b$ are locally defined functions on $M$
and $\tilde a$, $\tilde b$ are their lifts with respect to
$\sigma$, gives a well defined Poisson bracket on $\Cm$.
It is clear that equivalent formal deformations
have the same Poisson bracket.

\begin{definition} A deformation quantization (DQ) on
a symplectic manifold $(M,\omega)$ is a formal deformation of $\Cm$
whose Poisson bracket is equal to $\omega^{-1}$.
\end{definition}

\begin{definition} a) A star-product (SP) on $(M,\omega)$
is the structure of an associative algebra on the sheaf
$\CM$ with the multiplication of the form
\be{}
f\ast g=\sum_{i\geq 0} t^i\mu_i(f,g), \qquad f,g\in\Cm,
\ee{}
where all $\mu_i$, $i>0$, are bidifferential operators on $M$ vanishing
on constants, i.e. $\mu_i(f,g)=0$ if $f$ or $g$ is a constant,
$\mu_0(f,g)=fg$,
and $\mu_1(f,g)-\mu_1(g,f)=\{f,g\}$, the Poisson bracket inverse to
$\omega$.

b) Two star-products $(\CM,\mu')$ and $(\CM,\mu'')$
on $(M,\omega)$
are equivalent if there exists a power series
$B=1+tB_1+\cdots$, where
$B_i$ are differential operators vanishing on constants,
such that $\mu''(f,g)=B\mu'(B^{-1}f,B^{-1}g)$.
\end{definition}

It is clear that any SP $(\CM,\mu)$ defines
a DQ with the natural symbol map $\CM\to \Cm$,
$f_0+tf_1+\cdots\mapsto f_0$.

\begin{propn} The above assignment gives a one-to-one
correspondence between the equivalence classes of SP's and DQ's.
\end{propn}

\begin{proof}
Let $\bA_t$ be a DQ. Let us prove that it is equivalent to a
star-product. By definition of DQ, there exists an isomorphism of
$\C[[t]]$-modules $\CM\to\bA_t$ commuting with symbol maps. Let
$\mu=\mu_0+t\mu_1+\cdots$ be the multiplication in the sheaf $\CM$
being the pullback of the multiplication in $\bA_t$. In order to
prove that each $\mu_i$ is a bidifferential operator it is enough
to prove, according to the Peetre theorem, that
$supp(\mu(f,g))\subset supp(f)\cap supp(g)$ for any functions $f$
and $g$ on $M$. But this is obvious because if, for example, $f=0$
on an open set $U$, then $\mu(f,g)=0$ on $U$, since $\mu$ is a map
of sheaves. The same argument proves that two equivalent DQ's
correspond to equivalent SP's.
\end{proof}

\begin{example}\label{ex1} {\em Moyal-Weyl star-product.}
Let $U$ be an open set in a symplectic manifold $(M,\omega)$,
in which there exist Darboux coordinates,
$x_i$, $y_i$, $i=1,...,\half\dim M$, so that
$\omega=\sum_idy_i\wedge dx_i$.
Then,
$\omega^{-1}=\sum_i\partial/\partial y_i\wedge \partial/\partial x_i$
and for $f,g\in C^\infty_U$ the multiplication formula
\be{}
f\otimes g\mapsto m\exp\left(\frac{t}{2}\sum_i(\partial/\partial y_i
\ot\partial/\partial x_i-\partial/\partial x_i
\ot\partial/\partial y_i)\right)(f\otimes g),
\ee{}
where $m$ is the usual multiplication of functions,
defines a star-product on $(U, \omega)$.
This star-product is called the Moyal-Weyl star-product.
\end{example}

The following statement is well known and can be proven with the help
of Hochschild and de Rham cohomology.

\begin{propn} Locally, any star-product on $(M,\omega)$
is equivalent to the Moyal-Weyl star-product.
\end{propn}

\subsection{Polarized deformation quantization}

\begin{definition} a) A polarized deformation quantization (PDQ) on
a polarized symplectic manifold $(M,\omega,\PP)$
is a pair $(\bA_t,\bO_t)$ where $\bA_t$ is a DQ on $(M,\omega)$
and $\bO_t$ is a $t$-adically complete commutative subalgebra in $\bA_t$
such that $\sigma(\bO_t)=\OO_\PP$, functions constant along $\PP$.

b) Two PDQ's $(\bA_{t1},\bO_{t1})$ and $(\bA_{t2},\bO_{t2})$ of
$(M,\omega,\PP)$ are
equivalent if there exists an equivalence map of DQ's
$\bA_{t1}\to\bA_{t2}$ which takes $\bO_{t1}$ to $\bO_{t2}$.
\end{definition}

Since any DQ is equivalent to a SP, any PDQ is equivalent
to a triple $(\CM,\mu_t,\OO_t)$, where $(\CM,\mu_t)$ is
a SP and $\OO_t$ is a commutative subalgebra in
$(\CM,\mu_t)$.

\begin{propn}\label{propmaxq} Let $(\bA_t, \bO_t)$ be a PDQ. Then,

a) $\bO_t$ is a maximal commutative subalgebra in
$\bA_t$.

b) $\bO_t$ is a $t$-regular subalgebra
in $\bA_t$.
\end{propn}

\begin{proof}
a) Locally, the bracket $[a,b]=\tm(ab-ba)$, $a,b\in\bA_t$,
is a deformation of the Poisson bracket $\omega^{-1}$.
So, the statement easily follows from Proposition \ref{propmax}.

b) Follows from a).
\end{proof}

The following definition will play an auxiliary role
in the paper.

\begin{definition}
A {\em weakly polarized star-product} (wPSP) on
a polarized symplectic manifold $(M,\omega,\PP)$ is a triple
$(\CM,\mu_t,\OO_t)$, where $(\CM,\mu_t)$ is a SP,
$\OO_t$ is a $t$-adically complete $\Ct$-submodule in $\CM$ satisfying
the conditions:

a) $\sigma(\OO_t)=\OO_\PP$;

b) $\OO_t$ is a commutative subalgebra
in $\CM$ with respect to the usual multiplication in $\CM$;

c) $\mu_t$ being restricted to $\OO_t$ coincides with
the usual multiplication.
\end{definition}

\begin{definition}
We say that a wPSP $(\CM,\mu_t,\OO_t)$ is
a {\em polarized star-product} (PSP), if
for any $a\in\OO_t$, $b\in\CM$ the product
$\mu_t(a,b)$ coincides with the usual product in $\CM$.
\end{definition}

So, both wPSP and PSP are particular cases of a PDQ.
We are going to prove that, in fact, any PDQ is equivalent
to a PSP. But first we prove the following

\begin{lemma}\label{prop13}
Any PDQ is equivalent to a wPSP.
\end{lemma}

\begin{proof}
Since any deformation quantization is equivalent to a star-product,
it is enough to prove the following.
Let $\cO_t$ be a commutative subalgebra in a star-product
$(\CM,\mu_t)$ such that the triple
$(\CM,\mu,\cO_t)$ forms a PDQ, in particular,
$\mu_t(a,b)=\mu_t(b,a)$ for $a,b\in \cO_t$.
Then, there exists a differential operator
$D_t=1+tD_1+...$ on $M$ such that the multiplication
$\tilde\mu_t(f,g)=D_t^{-1}\mu(D_tf,D_tg)$ being restricted
to $D_t^{-1}\cO_t$ is the usual multiplication.

Suppose $\mu_t$ being restricted to $\cO_t$
coincides with the usual multiplication modulo $t^n$. Then,
$\mu(a,b)=ab+t^n\nu(a,b) \mod t^n$ for $a,b\in\OO_t$.
It follows from Proposition \ref{propmaxq} b) that
being considered modulo $t$, the bilinear form $\nu$
defines a Hochschild cocycle $\Bar{\nu}\in \Gamma(M,\D^2(\OP,\Cm))$
(see Proposition \ref{propKHR1}).

Since $\Bar{\nu}$ is commutative, it follows from Proposition \ref{propKHR1}
that it is a coboundary. So, there exists
a differential operator $\Bar{D}\in\Gamma(M,\D^1(\OP,\Cm))$
such that $d_{Hoch}\Bar{D}=\Bar{\nu}$.

Let $\tilde{D}$ be a lift of $\Bar{D}$ to a differential operator
on $M$. Let $D_n=1+t^n\tilde{D}$. It is easy to see
that $D_n^{-1}\cO_t$ is a commutative subalgebra in
the star-product $(\CM,\tilde{\mu}_t)$
with $\tilde{\mu}_t(a,b)=D_n^{-1}\mu_t(D_na,D_nb)$,
and $\tilde{\mu}_t(a,b)=ab$ modulo $t^{n+1}$ for $a,b\in D_n^{-1}\cO_t$.

By induction, we construct such differential operators
$D_n$ for all $n$.
Let $D=\Pi_{n=1}^\infty D_n$, $\cO'_t=D^{-1}\cO_t$, and
$\mu'_t(a,b)=D^{-1}\mu_t(Da,Db)$.
Then $D$ gives an isomorphism between PDQ's
$(\CM,\mu_t,\cO_t)$ and $(\CM,\mu'_t,\cO'_t)$,
and the second triple is a wPSP, which proves the proposition.
\end{proof}

\begin{propn}\label{prop13a}
For any wPSP $(\CM,\mu_t,\OO_t)$, there
exists a differential operator
$D=1+tD_1+\cdots$ on $M$ such that $Df=f$ for $f\in\OO_t$ and the
multiplication
$\mu'_t(a,b)=D^{-1}\mu_t(Da,Db)$ defines a PSP
$(\CM,\mu'_t,\OO_t)$.
\end{propn}
\begin{proof}
It is obvious that $\mu=\mu_t$ defines a PSP modulo $t$.
Proceeding by induction, we assume that there exists a
wPSP multiplication
$\mu'$ equivalent to $\mu$ and being a PSP modulo $t^n$
with respect to $\PP_t$. The proposition will be proved if
we find a differential operator, $D_n$, such that $D_n(f)=0$ for
all $f\in\OO_0$ and the multiplication
\be{}\label{form}
\mu''(a,b)=D^{-1}\mu'(Da,Db),
\ee{}
where $D=1+t^nD_n$,
defines a PSP modulo $t^{n+1}$
with respect to $\PP_t$.

Let
$$\mu'=\mu_0+t\mu'_1+\cdots +t^{n-1}\mu'_{n-1}+t^n\nu \quad \mod t^{n+1}.$$
By our assumption, elements $\mu'_1,...,\mu'_{n-1}$ are strongly polarized
and $\nu$ is polarized with respect to $\OO_t$ (see the definition before
Proposition \ref{propdop}).
It follows from associativity of $\mu'$ that
$$d_{Hoch}\nu(a,b,c)=\sum_{i+j=n}(\mu'_i(a,\mu'_j(b,c))-\mu'_i(\mu'_j(a,b),c)).$$
It is easy to check that each term in the right hand side is strongly
polarized, since all $\mu'_i$, $i=1,...,n-1$, are such.
So, $\nu$ satisfies the hypothesis of Proposition \ref{propdop}. Hence,
there exists
a polarized differential operator $D_n$ such that $\nu+d_{Hoch}D_n$
is a strongly polarized bidifferential operator. It is obvious
that the multiplication $\mu''$ defined in (\ref{form}) with $D_n$
just constructed is as required.
\end{proof}

\begin{cor}\label{cor1}
Any PDQ is equivalent to a PSP.
\end{cor}
\begin{proof}
Follows from Lemma \ref{prop13} and Proposition \ref{prop13a}.
\end{proof}

\begin{propn}\label{prop20}
a) Let $(\CM,\mu_t,\OO_t)$ be a PSP on $(M,\omega,\PP)$.
Let $\PP_t=(d\OO_t)^\perp$.
Then, $\PP_t$ is a deformation of $\PP$ and $\OO_t=\OO_{\PP_t}$.

b) Let $(\CM,\tilde\mu_t,\tilde\OO_t)$ be another PSP on
$(M,\omega,\PP)$ equivalent to
$(\CM,\mu_t,\OO_t)$ as a PDQ. Then, there exists a
formal automorphism of $M$ which takes $\PP_t$ to
$\tilde\PP_t=(d\tilde\OO_t)^\perp$.

c) Let $(\CM,\mu_t,\OO_t)$ and $(\CM,\tilde\mu_t,\OO_t)$
are two equivalent PSP's with the same $\OO_t$. Let $D=1+tD_1+\cdots$
gives an equivalence. Then,
there exists a decomposition, $D=D'e^{tX}$,
where $D'$ is a differential operator
identical on $\OO_t$ and $X$ is a formal vector field on $M$
taking $\OO_t$ to itself.
\end{propn}
\begin{proof}
a) follows from Proposition \ref{propdis}.

b). Let us put $X_0=0$. Then the automorphism $e^{tX_0}=Id$
takes $\PP_t$ to $\tilde\PP_t$
modulo $t$. Suppose we have constructed a formal vector field $X_k$ such that the formal automorphism
$e^{tX_{k-1}}$ of $M$ takes $\PP_t$ to $\tilde\PP_t$ modulo $t^{k}$. Then, replacing $\PP_t$ by
$e^{tX_{k-1}}\PP_t$ we may assume that $\PP_t$ and $\tilde\PP_t$ coincide
modulo $t^k$. The proposition will
be proved, if we show that it is possible to find a vector field
$Y$ such that $e^{t^kY}$ takes $\PP_t$ to
$\tilde\PP_t$ modulo $t^{k+1}$.

Since, by our assumption, the SP's $(\CM,\mu_t)$ and
$(\CM,\tilde\mu_t)$ are equivalent and coincide modulo
$t^k$, there exists a differential operator $1+t^kD_k$ realizing
that equivalence.
Since $1+t^kD_k$ takes $\OO_t$ to $\tilde\OO_t$ and on the both
of these subalgebras the respecting
multiplications $\mu_t$ and $\tilde\mu_t$ are trivial,
$D_k$ being restricted to $\OO_0$ is a derivation
from $\OO_0$ to $\Cm$. Let $Y$ be an extension of that restricted
$D_k$ to a derivation on $\Cm$. It is clear that $e^{t^kY}$
takes $\PP_t$ to $\tilde\PP_t$ modulo $t^{k+1}$.

c). The operator $D$
being restricted to $\OO_t$ is a formal automorphism of $\OO_t$.
Since formal automorphisms form a
pro-unipotent group, there exists $X'\in Der(\OO_t)$
such that $D$ being restricted to $\OO_t$ is equal to
$e^{tX'}$. Let $X\in\Der(\Cm)$ be a lift of $X'$.
We put $D'=De^{-tX}$ which is obviously identical on $\OO_t$.
\end{proof}

\begin{example}\label{examp2}  {\em Moyal-Wick PSP.}
Let $(M,\omega,\PP)$ be a polarized symplectic manifold.
Let $U$ be an open set in $M$
where there exist Darboux coordinates,
$x_i$, $y_i$, $(dx_i)^\perp=\PP$, $i=1,...,\half\dim M$, so that
$\omega=\sum_idy_i\wedge dx_i$.
Then,
$\omega^{-1}=\sum_i\partial/\partial y_i\wedge \partial/\partial x_i$
and for $f,g\in C^\infty_U$ the multiplication formula
\be{}
f\otimes g\mapsto m\exp\left(t\sum_i\partial/\partial y_i
\ot\partial/\partial x_i\right)(f\otimes g),
\ee{}
where $m$ is the usual multiplication of functions,
defines a PSP on $(U, \omega, \PP)$.
This PSP is called the Moyal-Wick polarized star-product.

Note that the functions $a_i,f_i$ satisfying the Darboux
relations with respect to the Poisson bracket $\omega^{-1}$
also satisfy the Darboux relations with respect to
the deformed bracket $\tm[\cdot,\cdot]$, where
$[\cdot,\cdot]$ is the commutator of the Moyal-Wick PSP.

Let us remark that the Moyal-Weyl SP from Example \ref{ex1}
constructed using the same Darboux
coordinates
gives just a wPSP but not a PSP.
\end{example}

\begin{propn} Locally, any PSP on $(M,\omega,\PP)$ is
equivalent to the Moyal-Wick PSP.
\end{propn}

\begin{proof}
Let us prove that any two PSP's are locally equivalent.
Since any SP's are locally equivalent, we may suppose
that there are given two PSP's, $(\CM,\mu_t,\OO_t)$ and
$(\CM,\mu_t,\tilde\OO_t)$, with the same
multiplication and different polarizations,
and we have to prove that they are, locally, equivalent.
Let $x_i$, $y_i$ are Darboux coordinates with respect to
$\omega$ such
that $(dx_i)^\perp=\PP$. By Proposition \ref{prop11}, there exist
their lifts $x_{it}$,
$y_{it}$ and $x'_{it}$, $y'_{it}$
which satisfy the Darboux relations with respect to the bracket
$[a,b]=\frac{1}{t}(\mu_t(a,b)-\mu_t(b,a))$, and
$x_{it}\in\OO_t$, $x'_{it}\in\tilde\OO_t$. By Proposition \ref{prop12},
there exists, locally, an inner
automorphism of the SP $(\CM,\mu_t)$ that takes
$x_{it}$, $y_{it}$ to $x'_{it}$, $y'_{it}$. It follows that
this automorphism takes $\OO_t$ to $\tilde\OO_t$.
\end{proof}

\section{Characteristic classes of PDQ's and PSP's}

\subsection{Extension class associated with a PDQ}

Let $(\bA_t,\bO_t)$ be a PDQ on a polarized symplectic manifold
$(M,\omega,\PP)$.
Since any PDQ is equivalent to
a PSP, the sheaf $\bO_t$ is isomorphic to
$\OO_{\PP_t}$ for some deformed distribution $\PP_t$.
Thus, the sheaves  $Der(\bO_t)$ and
$\Omega^{1cl}_{\bO_t}$ are well defined
(see Section \ref{subsec2.2}).
Let
\be{}\label{exten}
F(\bA_t,\bO_t)=\{b\in\bA_t; [b,\bO_t]\subset\bO_t\},
\ee{}
where $[\cdot,\cdot]$ denotes the commutator in $\bA_t$.
It is clear that $F(\bA_t,\bO_t)$ is a sheaf of Lie algebras
with the bracket $\ft[\cdot,\cdot]$ and the center $\bO_t$.
Moreover, any element $b\in F(\bA_t,\bO_t)$
determines the derivation $\ft[b,\cdot]$
of $\bO_t$ and, due to Proposition \ref{prop11},
this correspondence defines an epimorphism
$\sigma:F(\bA_t,\bO_t)\to Der(\bO_t)$.

We consider $F(\bA_t,\bO_t)$ as a left $\bO_t$-module
with respect to multiplication in $\bA_t$. As a Lie algebra
sheaf, $F(\bA_t,\bO_t)$ is an extension of $Der(\bO_t)$.

So, we have the following exact sequence of Lie algebras and
$\OO_t$-modules:
\be{}\label{defoext}
\begin{CD}
0 @>>> \bO_t @>\iota>> F(\bA_t,\bO_t) @>\sigma>> Der(\bO_t) @>>> 0.
\end{CD}
\ee{}
According to the terminology of \cite{BB}, \cite{BK},
$F(\bA_t,\bO_t)$ is called {\em a $\bO_t$-extension of $Der(\bO_t)$}.

We say that a map of Lie algebras and $\bO_t$-modules,
$s: Der(\bO_t)\to F(\bA_t)$, given over an open set of $M$
is a splitting of (\ref{defoext}), if $s\sigma=id$.
Since $(\bA_t,\bO_t)$ can be realized as a PSP and, locally,
there exist Darboux coordinates with respect to $\ft[\cdot,\cdot]$,
the sequence (\ref{defoext}) locally splits
(see the next subsection, where splittings are presented
with the help of Darboux coordinates).

\begin{lemma}\label{splitting}
Let $s$ and $s'$ are two splittings of (\ref{defoext}) over an
open set of $M$. Then, $s-s'\in \Omega^{1cl}_{\bO_t}$.
\end{lemma}
\begin{proof} Direct calculation.
\end{proof}

Let us define the extension class of (\ref{exten})
in the following way.
Let $\{U_\alpha\}$ be an open covering of $M$
such that over each $U_\alpha$ there is a splitting, $s_\alpha$,
of $F(\bA_t,\bO_t)$. By Lemma \ref{splitting},
$f_{\alpha,\beta}=s_\beta-s_\alpha$
is a section of $\Omega^{1cl}_{\bO_t}$ over $U_\alpha\cap U_\beta$.
We define $\cle(F(\bA_t,\bO_t))$ as the element of
$H^1(M,\Omega^{1cl}_{\bO_t})$
represented by the collection $\{f_{\alpha,\beta}\}$.
One can prove that given $\bO_t$, the extension class $\cle(F(\bA_t,\bO_t))$
determines a $\bO_t$-extension of $Der(\bO_t)$ up to
equivalence.

We will denote the element
$\cle(F(\bA_t,\bO_t))$
by $\cle(\bA_t,\bO_t)$ and call it
the {\em extension class} of the PDQ $(\bA_t,\bO_t)$.

In the next subsection, the extension class
of a PSP, $(\CM,\mu_t,\OO_t)$, will be represented
as an element of $\Hwt$, $\PP_t=\PP_{\OO_t}$,
with the help of a characteristic 2-form associated with that PSP.

\subsection{Characteristic 2-form associated with a PSP}
\label{subsecextcl}

Let $(\CM,\mu_t,\OO_t)$ be a PSP
that we denote for shortness by $(\mu_t,\OO_t)$.
We denote $\PP_t=\PP_{\OO_t}$.
Then the extension $F(\mu_t,\OO_t)$
coincides as a left $\OO_t$-module
with a $\OO_t$-submodule of $\CM$ with respect to
the usual multiplication in $\CM$.
In this case, the local splittings of $F(\mu_t,\OO_t)$
are differential forms of
$Hom_{\OO_t}(Der(\OO_t),\CM)=\PP^\perp$.
These forms can be described
explicitly by Darboux coordinates.

Let $\{U_\alpha\}$ be an open covering of $M$ such that
each $U_\alpha$ has Darboux coordinates
$x_{\alpha i}$, $y_{\alpha i}$, $x_{\alpha i}\in \OO_t$,
$i=1,...,n$,
with respect to
the bracket $[a,b]=\ft(\mu_t(a,b)-\mu_t(b,a))$
(in particular, $[y_{\alpha i},x_{\alpha j}]=\delta_{ij}$).
By Proposition \ref{prop11} such a covering exists.
Then, on $U_\alpha$, the $\OO_t$-submodule $F(\mu_t,\OO_t)\subset\CM$
is equal to
$\OO_t\oplus\OO_ty_{\alpha 1}\oplus\cdots\oplus\OO_ty_{\alpha n}$.

Splittings $s_\alpha$ may be taken by the condition
$x_{\alpha i}\mapsto y_{\alpha i}$, and the
corresponding forms are
\be{}\label{Darsplit}
s_\alpha=\sum_i y_{\alpha i}dx_{\alpha i}.
\ee{}
By Lemma \ref{splitting}, $ds_\alpha=ds_\beta$
on $U_\alpha\cap U_\beta$.
So, forms $s_\alpha$ define the global 2-form
\be{}\label{gform}
\omega_t\in \omega_0+t\Gamma(M,d\PP_t^\perp), \quad\quad
\omega_t=ds_\alpha=\sum_i dy_{\alpha i}\wedge dx_{\alpha i}.
\ee{}
\begin{lemma}\label{lemrepcl}
This form
represents the extension class $\cle(\mu_t,\OO_t)$
by the isomorphism (\ref{fundiso1}).
\end{lemma}
\begin{proof} Clear.
\end{proof}
Lemma \ref{splitting} shows that if we take other splittings
of $F(\mu_t,\OO_t)$,
the procedure above gives the same form $\omega_t$.
So, in the case of PSP, it is well defined
not only the class
$\cle(\mu_t,\OO_t)$ but also its
representative in
$\omega_0+t\Gamma(M,d\PP_t^\perp)$ that we call the
{\em characteristic 2-form} of the PSP and denote it by
$\clP(\mu_t,\OO_t)$.

Given a SP, $(\CM,\mu_t)$, let us denote
$[a,b]_{\mu_t}=\frac{1}{t}(\mu_t(a,b)-\mu_t(b,a))$, $a,b\in\CM$.
The bracket $[,]_{\mu_t}$ is a deformation of the initial
Poisson bracket on $M$.

\begin{propn}\label{propcom}
Let $(\mu_t,\OO_t)$, $(\tilde\mu_t,\OO_t)$ be PSP's
with the same $\OO_t$.
Then, $\clP(\mu_t,\OO_t)=\clP(\tilde\mu_t,\OO_t)$ if and only if
$[,]_{\mu_t}=[,]_{\tilde\mu_t}$.
\end{propn}
\begin{proof}
If $\clP(\mu_t,\OO_t)=\clP(\tilde\mu_t,\OO_t)$,
then forms (\ref{Darsplit})
can be taken to be the same. This implies that
there exist common Darboux coordinates with respect to
$[,]_{\mu_t}$ and $[,]_{\tilde\mu_t}$.
\end{proof}

\begin{propn}\label{propeq}
Let $(\mu_t,\OO_t)$, $(\tilde\mu_t,\OO_t)$ be two PSP's
with the same $\OO_t$. Let $\PP_t=\PP_{\OO_t}$ and
$\Gamma\PP_t$ denote the module of global sections of $\PP_t$.
Then, the 2-forms $\clP(\mu_t,\OO_t)$, $\clP(\tilde\mu_t,\OO_t)$
are lying on the same orbit of $e^{t\Gamma\PP_t}$ if and only if
there exists a formal differential operator $D=1+tD_1+\cdots$
identical on $\OO_t$ such that
$\tilde\mu_t(a,b)=D^{-1}\mu_t(Da,Db)$.
\end{propn}

\begin{proof}
Let us denote $\omega_t=\clP(\mu_t,\OO_t)$,
$\tilde\omega_t=\clP(\tilde\mu_t,\OO_t)$.
Assume such $D$ exists. Let $s_\alpha$ be forms
from (\ref{Darsplit}) corresponding to $(\mu_t,\OO_t)$.
Let $\tilde s_\alpha$ be forms obtained
from $s_\alpha$ by applying $D$.
Lemmas \ref{splitting} and \ref{lemma2.2} imply
that there exist $f_{\alpha,\beta}\in\OO_t$
over $U_\alpha\cap U_\beta$ such that
$s_\alpha-s_\beta=df_{\alpha,\beta}$. Since $D$ acts on $\OO_t$
trivially, one has
$\tilde s_\alpha-\tilde s_\beta=df_{\alpha,\beta}$, too.
This means that $\tilde s_\alpha-s_\alpha$ does not depend
on $\alpha$ and give a global
form, $b$, of $\PP_t^\perp$. Since $\omega_t=ds_\alpha$,
$\tilde\omega_t=d\tilde s_\alpha$, one has
$\tilde\omega_t=\omega_t+db$.
By Lemma \ref{lemma3}, there exists a formal
automorphism, $e^{tY}$, $Y\in\Gamma\PP_t$, taking
$\omega_t$ to $\tilde\omega_t$.

Conversely, let us suppose that $\omega_t$ and $\tilde\omega_t$
are on the same orbit of $e^{t\Gamma\PP_t}$.
Assume that we have found a differential operator
identical on $\OO_t$ which transforms $\tilde\mu_t$ to
a multiplication $\tilde\mu'_t$ that is equal to $\mu_t$
modulo $t^k$, i.e.
\be{}\label{modtk}
\tilde\mu'_t-\mu_t=t^k\nu+\cdots.
\ee{}
By the previous part of the proof, the corresponding form
$\tilde\omega'_t$ is also lying on the same orbit as $\omega_t$.
The proposition will be proved if we find a differential operator
$1+t^kD_k$, $D_k(\OO_t)=0$, which transforms
$\tilde\mu'_t$ to $\mu_t$ modulo $t^{k+1}$.
Let us prove that.

Since $\tilde\mu'_t=\mu_t$ mod $t^k$, we can choose, by
Proposition \ref{prop11a}, systems of Darboux coordinates with
respect to $[,]_{\tilde\mu'_t}$ and $[,]_{\mu_t}$ that
coincides modulo $t^k$.
It follows that
$\tilde\omega'_t=\omega_t$ mod $t^k$. Hence,
there is a formal automorphism $e^{t^kX}$, $X\in\PP_t$,
which takes $\tilde\omega'_t$ to $\omega_t$.
Applying $e^{t^kX}$ to $\tilde\mu'_t$, we obtain a multiplication
that is still equal to $\mu_t$ modulo $t^k$ but
whose characteristic form is equal to $\omega_t$.
So, we come to
the situation when the multiplications $\tilde\mu'_t$
and $\mu_t$ have the same characteristic form $\omega_t$.

By Proposition \ref{propcom},
$\tilde\mu'_t$ and $\mu_t$ have the same commutator.
This implies that
bidifferential operator $\nu$ in
(\ref{modtk}) is commutative. Moreover, it is a Hochschild
strongly polarized cocycle (see Section 2). So,
there exists a polarized differential operator $D_k$
such that $d_{Hoch}D=\nu$.
It follows that transformation $1+t^kD_k$ applying
to $\tilde\mu'_t$ gives a multiplication equal to
$\mu_t$ modulo $t^{k+1}$ and identical on $\OO_t$.
\end{proof}

\subsection{Characteristic pairs for PSP's and PDQ's}\label{subsclass}

Let $(M,\omega,\PP)$ be a polarized symplectic manifold.
Let us denote by $Aut(M)$ the group of formal automorphisms
of $M$ and by $\YY=\YY(M,\omega,\PP)$ the set of pairs
$(\omega_t,\PP_t)$, where $\omega_t=\omega_0+t\omega_1+\cdots$
is a formal symplectic form being a deformation of $\omega=\omega_0$
and $\PP_t$ is a polarization of $\omega_t$
being a deformation of $\PP$. The group $Aut(M)$ naturally
acts on $\YY$.

It is natural to assign
to a PSP $(\mu_t,\OO_t)$ on $(M,\omega,\PP)$
a pair $(\omega_t,\PP_t)\in\YY$,
where $\omega_t=\clP(\mu_t,\OO_t)$ and $\PP_t=\PP_{\OO_t}$.
So, we obtain the map
\be{*}
\tau:\{\mbox{PSP's}\}\to\YY.
\ee{*}
\begin{propn}\label{eqcl}
Two PSP's $(\mu_t,\OO_t)$ and $(\tilde\mu_t,\tilde\OO_t)$ are
equivalent if and only if $\tau(\mu_t,\OO_t)$ and
$\tau(\tilde\mu_t,\tilde\OO_t)$ are lying on the same orbit
in $\YY$ under the $Aut(M)$-action.
\end{propn}
\begin{proof}
Let $(\mu_t,\OO_t)$ and $(\tilde\mu_t,\tilde\OO_t)$ be
equivalent. Let us prove that the pairs
$\tau(\mu_t,\OO_t)$ and
$\tau(\tilde\mu_t,\tilde\OO_t)$ are lying on the same orbit.
By Proposition \ref{prop20} b), c)
one can find a formal automorphism of $M$
such that after its applying we come to the situation when
$(\tilde \mu_t,\tilde \OO_t)$ turns into a PSP,
$(\tilde \mu_t,\OO_t)$, with the same $\OO_t$ as in
$(\mu_t,\OO_t)$ and the equivalence
morphism from $(\mu_t,\OO_t)$ to
$(\tilde \mu_t,\OO_t)$ is given by a differential
operator identical on $\OO_t$. Now the statement follows from
Proposition \ref{propeq}.

Conversely, suppose that for
$(\mu_t,\OO_t)$ and $(\tilde\mu_t,\tilde\OO_t)$
the corresponding pairs
$(\omega_t,\PP_t)$, $(\tilde\omega_t,\tilde\PP_t)$
lie on the same orbit.
Let us prove that those PSP's are equivalent.
Let $B$ be a formal automorphism of $M$ sending
$\tilde\PP_t$ to $\PP_t$.
Applying $B$ to $(\tilde\mu_t,\tilde\OO_t)$,
we come to the case when $\tilde\PP_t=\PP_t$.
So, we may suppose that
$(\tilde\omega_t,\tilde\PP_t)=(\tilde\omega_t,\PP_t)$.
Since the pairs $(\omega_t,\PP_t)$, $(\tilde\omega_t,\PP_t)$
lie on the same orbit,
there exists $X\in\Gamma\PP_t$ such that
$\tilde\omega_t=e^{tX}\omega_t$. Now the statement
follows from Proposition \ref{propeq}.
\end{proof}

Let us denote by $\bigl[\YY\bigr]$ the set of orbits in
$\YY=\YY(M,\omega,\PP)$.

\begin{cor}\label{corin}
The map $\tau$ induces the embedding
\be{*}
\clPQ: \{\mbox{classes of PDQ's}\}\to \bigl[\YY\bigr].
\ee{*}
\end{cor}
\begin{proof}
Let $(\bA_t,\bO_t)$ be a PDQ on $(M,\omega,\PP)$.
Then, by Corollary \ref{cor1}, there exists a PSP,
$(\mu_t,\OO_t)$, equivalent to $(\bA_t,\bO_t)$.
We put $\clPQ(\bA_t,\bO_t)=[\tau(\mu_t,\OO_t)]$, the orbit
passing through the pair $\tau(\mu_t,\OO_t)$.
By Propositions \ref{prop20} and \ref{eqcl}, this map is
correctly defined, i.e. does not depends on the
choice of an equivalent PSP. Proposition \ref{eqcl}
also shows that being descended to equivalence
classes of PDQ's $\clPQ$ is embedding.
\end{proof}

In the next section we will prove that any element
of $\YY$ is equal to $\tau(\mu_t,\OO_t)$ for a PSP $(\mu_t,\OO_t)$,
which implies
that the map $\clPQ$ is, in fact, an isomorphism.

\section{Existence of polarized deformation quantizations
and relation between the extension and Fedosov classes of a PDQ}

Let $(M,\omega_0)$ be a symplectic manifold.
It is known that all equivalence classes of deformation quantizations
(DQ) on
$M$ with the Poisson bracket $\omega_0^{-1}$
can be obtained by the Fedosov method. According to this method,
starting with a symplectic connection,
one constructs a flat connection, $D$, (called the Fedosov connection)
in the Weyl algebra defined on
the cotangent bundle to $M$ via the Poisson
bracket $\omega_0^{-1}$.
The quantized algebra, $\bA_t$, is realized as the subalgebra of
flat sections in the Weyl algebra. The Weyl curvature of $D$ (see below),
being a closed scalar 2-form of the form
$\theta_t=\omega_0+t\omega_1+\cdots$,
defines the Fedosov class
\be{}\label{formtheta}
cl_F(\bA_t)=[\theta_t]\in [\omega_0]+tH^2(M,\C[[t]]).
\ee{}
It is also known that the correspondence $\bA_t\mapsto cl_F(\bA_t)$
is a bijection between the set of equivalence classes of DQ's
on $(M,\omega_0)$ and the set $[\omega_0]+tH^2(M,\C[[t]])$,
\cite{Fe}, \cite{NT}, \cite{Xu}.

Let $(M,\omega,\PP)$ be a polarized symplectic manifold
and $(\omega_t,\PP_t)\in\YY(M,\omega,\PP)$ a deformation
of the pair $(\omega,\PP)$, see Subsection \ref{subsclass}.
We adapt the Fedosov method to construct a PSP,
$(\mu_t,\OO_t)$, such that $\tau(\mu_t,\OO_t)=(\omega_t,\PP_t)$.
We start with a $\PP_t$-symplectic connection, $\nabla$,
and construct the Fedosov connection on the same Weyl
algebra.
We will see that
by realizing the Fedosov scheme in presence of a polarization,
the form $\omega_t$ appears as
a so-called Wick curvature of the Fedosov connection.
Moreover, $\omega_t$
differs from the Weyl curvature of that Fedosov connection
by the form $\ftt tr(\nabla^2|_{\PP_t})$.

\subsection{Some notations}

Let $\EE$ be a formal vector bundle over $M$, i.e. a
free $\CM$-module of finite rank over $M$.
Denote by $T^k(\EE)$ the $k$-th tensor power of $\EE$ over $\CM$ and by
$T(\EE)$ the corresponding tensor algebra completed in the
$\{\EE,t\}$-adic topology.
Similarly we define the completed symmetric algebra
$S(\EE)$. For  a subbundle $\PP$ of $\EE$, we denote by $sym_\PP:
S(\PP)\to T(\EE)$ the natural embedding of $\CM$-modules defined by
symmetrization.

Let $\La(\EE)$ be the exterior algebra of $\EE$ over $\CM$. We will
consider $T(\EE)\ot\La(\EE)$ as a graded super-algebra
regarding a section $x\in T(\EE)\ot\La^k(\EE)$ of degree $k$ even (odd)
if $k$ is even (odd).

Denote by $\delta_{T(\EE)}$ the continuous $\CM$-linear
derivation of $T(\EE)\ot\La(\EE)$ defined by the map $T^1(\EE)\ot 1\to
1\ot\La^1(\EE)$, $v\ot 1\mapsto 1\ot v$, $v\in \EE$ is a section. It is
clear that $\delta_{T(\EE)}$ is a derivation of degree $1$ and
$\delta_{T(\EE)}^2=0$.
It is easy to see that for
any subbundle $\PP\subset \EE$, the map $\delta_{T(\EE)}$ being
restricted to $S(\PP)\ot\La(\PP)$ via the embedding $sym_\PP\ot
id_\PP$ gives a derivation of the
algebra $S(\PP)\ot\La(\PP)$; we denote it by $\delta_\PP$.

On the algebra $S(\PP)\ot\La(\PP)$, there is another derivation,
$\delta^*_\PP$, of degree $-1$ generated by the map
$1\ot\La^1(\PP)\to S^1(\PP)\ot1$,
$1\ot v\to v\ot 1$, $v\in\PP$.
It is easy to check that $(\delta^*_\PP)^2=0$ and
$[\delta_\PP,\delta^*_\PP]=\delta_\PP\delta^*_\PP+\delta^*_\PP\delta_\PP=deg$,
where $deg$ is
the derivation assigning to an element $x\in S^p(\PP)\ot\La^q(\PP)$
the element $(p+q)x$.

Let $\EE$ be presented as a direct sum of $\CM$-submodules,
$\EE=\PP\oplus\QQ$.
It is obvious that the derivations $\delta_\PP$, $\delta^*_\PP$,
$\delta_\QQ$, $\delta^*_\QQ$ induce derivations on
the algebra $\Sbf(\PP,\QQ)=(S(\PP)\ot S(\QQ))\ot (\La(\PP)\ot\La(\QQ))$
that we will identify in the natural way with
the algebra $(S(\PP)\ot S(\QQ))\ot \La(\EE)$.
We put $\delta_{\PP,\QQ}=\delta_\PP+\delta_\QQ$
and $\delta^*_{\PP,\QQ}=\delta^*_\PP+\delta^*_\QQ$.

Let us define the operator $\delta^{-1}_{\PP,\QQ}$ on
$\Sbf(\PP,\QQ)$ in the following way. We put $\delta^{-1}_{\PP,\QQ}(x)=0$
for $x\in \CM$ and
$\delta^{-1}_{\PP,\QQ}(x)=(1/(p+r+q)\delta_{\PP,\QQ}^*(x)$ for
$x\in(S^p(\PP)\ot S^r(\QQ))\ot \La^q(\EE)$, $p+r+q>0$.
There is the obvious relation
\be{}\label{relde}
\delta_{\PP,\QQ}\delta_{\PP,\QQ}^{-1}+
\delta_{\PP,\QQ}^{-1}\delta_{\PP,\QQ}= \mbox{ projection on $\Sbf^+(\PP,\QQ)$
along $\CM$},
\ee{}
where $\Sbf^+(\PP,\QQ)$ is the closure of $\oplus_{p+r+q>0}(S^p(\PP)\ot
S^r(\QQ))\ot \La^q(\EE)$.

One has the embedding
\be{}\label{rels}
sym_\PP\ot sym_\QQ\ot id:
(S(\PP)\ot S(\QQ))\ot \La(\EE)\to T(\EE)\ot\La(\EE).
\ee{}
It is clear that
$\delta_{\PP,\QQ}$ coincides with
the restriction of $\delta_{T(\EE)}$ to $\Sbf(\PP,\QQ)$ via this
embedding.

\subsection{The Fedosov algebra}
Let $\ff:\EE\ot\EE\to\CM$ be a $\CM$-linear skew-symmetric nondegenerate
form and $I$
the closed ideal in $T(\EE)$ generated by relations
\be{}\label{relI}
x\ot y-y\ot x=t\ff(x,y).
\ee{}
We call $\WW(\EE)=T(\EE)/I$ {\em the Weyl
algebra} and $\Wb=\Wb(\EE)=\WW\ot\La(\EE)$ {\em the Fedosov algebra}
over $\EE$. The derivation $\delta_{T(\EE)}$ on $T(\EE)\ot\La(\EE)$
induces a derivation of
$\Wb$. Indeed, $\delta_{T(\EE)}$ applied to the both
sides of (\ref{relI}) gives zero. We denote this derivation
by $\delta$.

Let $\EE=\PP\oplus\QQ$ be a decomposition into $\CM$-modules.

Define the Wick map, $\wi_{\PP,\QQ}$, as the composition
$\Sbf(\PP,\QQ)\to T(\EE)\ot\La(\EE)\to \Wb$, where the first map is
(\ref{rels}) and the second one is the projection.
By the PBW theorem $\wi_{\PP,\QQ}$
is an isomorphism of $\CM$-modules.

Due to the isomorphism $\wi_{\PP,\QQ}$, the operators
$\delta_{\PP,\QQ}$
and $\delta^{-1}_{\PP,\QQ}$ are carried over from $\Sbf(\PP,\QQ)$ to $\Wb$.
We retain for them the same notation. Note that while
$\delta_{\PP,\QQ}$ does not depend on the decomposition of $\EE$ and
coincides with the derivation $\delta$
induced from $T(\EE)\ot\La(\EE)$, the operator
$\delta^{-1}_{\PP,\QQ}$ is not a derivation and does depend on the
decomposition. In particular, one can suppose that the decomposition is
trivial, $\EE=\EE\oplus 0$. In this case we denote
$\delta^{-1}_\EE=\delta^{-1}_{\EE,0}$.

\begin{propn}\label{prop1.1} One has
$$H(\Wb,\delta)=\CM.$$ Moreover, if $x\in \WW(\EE)\ot\La^{k>0}(\EE)$
then $y=\delta^{-1}_{\PP,\QQ}x$ is such that $\delta y=x$ for any
decomposition $\EE=\PP\oplus\QQ$.
\end{propn}
\begin{proof} Follows from (\ref{relde}).
\end{proof}

\subsection{Lie subalgebras in $\WW$}

Let $\EE=\PP\oplus\QQ$ be a decomposition.
We say that $x\in \Wb$ has $\wi_{\PP,\QQ}$-degree $(p,q)$ if
$\wi_{\PP,\QQ}^{-1}(x)\in (S^p(\PP\ot S^q(\QQ))\ot\La(\EE)$.
We say that $x\in \Wb$ has $\wi_{\PP,\QQ}$-degree $k$ if
$\wi_{\PP,\QQ}^{-1}(x)\in \oplus_{p+q=k}(S^p(\PP)\ot S^q(\QQ))\ot\La(\EE)$.
We define the $\wi_\EE$-degree as the $\wi_{\EE,0}$-degree for the trivial
decomposition $\EE=\EE\oplus 0$.

Let $\g$ be a sheaf of Lie algebras acting on $\EE$.
We call a $\CM$-linear map $\lambda:\g\to\WW$ {\em a realization} of
$\g$, if it is a Lie algebra morphism ($\WW$ is considered as a Lie algebra with respect
to the commutator $\ft [\cdot,\cdot]$)
and for any $x\in\g$ and $v\in\EE$
one has $x(v)=\ft[\lambda(x),v]$. It is easy to check that
any two realizations differ
by a Lie algebra morphism of $\g$ to the center of $\WW$, so
if $\g$ is a sheaf of semisimple Lie algebras,
there is not more than one realization of $\g$.

Denote by ${\frak{sp}}(\EE)$ the sheaf of symplectic Lie algebras
with respect to $\ff$. Since ${\frak{sp}}(\EE)$ is semisimple, there is a unique
realization  $\rho_\EE:{\frak{sp}}(\EE) \to \WW$. The image of this realization
consists of elements having $\wi_\EE$-degree two.

Let $\EE=\PP\oplus\QQ$ be a decomposition into Lagrangian subsheaves.
Denote by ${\frak{sp}}(\PP,\EE)$ the subsheaf of ${\frak{sp}}(\EE)$
preserving $\PP$.
It is easy to check that ${\frak{sp}}(\PP,\EE)$ can be realized as the subset
of elements of $\WW$ having $\wi_{\PP,\QQ}$-degree $(1,1)$ and $(2,0)$.
Denote this realization by $\rho_{\PP,\EE}:{\mathfrak{sp}}(\PP,\EE)\to \WW$.
On the other hand, ${\mathfrak{sp}}(\PP,\EE)$ is realized in
$\WW$ by $\rho_\EE$.

\begin{lemma}\label{lem1.1}
Let $b\in {\frak{sp}}(\PP,\EE)$.
Then
$$\rho_\EE(b)=\rho_{\PP,\EE}(b)+\frac{t}{2}trace(\Bar b),$$
where $\Bar b$ is $b$ restricted to $\PP$.
\end{lemma}

\begin{proof} Direct computation using the fact that
$\rho_{\PP,\EE}(\Bar b)$ is
$(1,1)$ $\wi_{\PP,\QQ}$-degree component of $\rho_{\PP,\EE}(b)$ in any decomposition
$\EE=\PP\oplus\QQ$.
\end{proof}

\subsection{Filtrations on $\WW$}

We define two decreasing filtrations on $\WW$
numbered by nonnegative integers.

The $T$-filtration $F^T_\bullet\WW$ is defined as follows.
We ascribe to the elements of $\EE$ degree 1 and to $t$ degree 2.
Then $F^T_n\WW$  consists of elements of $\WW$
having the leading term of total degree $\geq n$.

The $\PP$-filtration, $F^\PP_\bullet\WW$, is firstly defined on $S(\PP)\ot S(\QQ)$
by the subsets
$F^\PP_n=\oplus_{i\geq n}S^i(\PP)\ot S(\QQ)$, $n=0,1,...$,
and carried over to $\WW$ via the Wick
isomorphism.

We extend those filtrations to $\Wb$ in the natural way standing,
for example, $F^T_n\Wb=F^T_n\WW\ot\La(\EE)$.
We will use the following mnemonic notation. To point out, for example,
that a section $x\in \Wb$ belongs to $F^T_n\Wb$ we write
$F^T(x)\ge n$.

In the following we denote
$\Sbf(\PP)=S(\PP)\ot\La(\EE)$ embedded in $\Sbf(\PP,\QQ)$
as $(S(\PP)\ot 1)\ot\La(\EE)$.

\begin{propn}\label{prop1.2}
Let $\EE=\PP\oplus \QQ$ be a decomposition of $\EE$ into
Lagrangian subsheaves. Then

a) The Wick map $\wi_{\PP,\QQ}:\Sbf(\PP,\QQ)\to\Wb$
has the following property:
for $a\in \Sbf(\PP)$
and arbitrary $c\in\Sbf(\PP,\QQ)$ one has
$ac=\wi(a)\wi(c)$.

 The filtrations on $\Wb$ have the properties:

b) for $x,y\in \Wb$, if $F^\PP(x)\ge k$, then $F^\PP(xy)\ge k$;

c) $F^\PP(\delta^{-1}_{\PP,\QQ}x)\geq F^\PP(x)$;

d) $F^T(\delta^{-1}_{\PP,\QQ}x)>F^T(x)$.
\end{propn}

\begin{proof} Clear.
\end{proof}

\subsection{Fedosov's construction in the Wick case}
\label{FC}

Let $(M,\omega_t,\PP_t)$ be a formal polarized symplectic manifold.
We will write for shortness $\omega=\omega_t$, $\PP=\PP_t$.
Let us denote $\TT=\TT^{\C}_M[[t]]$.
It is easy to prove that there exists on $M$ a
Lagrangian subbundle
$\QQ\subset\TT$ such that $ \TT=\PP\oplus\QQ$.

In the following we set $\EE=\TT^*$ and consider
the Fedosov algebra $\Wb=\Wb(\EE)$ with respect
to $\ff$ being the Poisson bracket inverse to $\omega$.
The decomposition $\EE=\PP^\perp\oplus\QQ^\perp$ is Lagrangian
with respect to this $\ff$.

Let $\nabla$ be a $\PP$-symplectic connection on $M$
(see Section \ref{PSC}).
Then the induced connection $\nabla:\EE\to\EE\ot\La^1\EE$
on $\EE$ preserves $\PP^\perp$, i.e.
$\nabla(\PP^\perp)\subset\PP^\perp\ot\La^1(\EE)$,
and is flat on $\PP^\perp$ along $\PP$, i.e.
for any $X,Y\in\PP$ one has
$(\nabla_X\nabla_Y-\nabla_Y\nabla_X-\nabla_{[X,Y]})(\PP^\perp)=0$.

We will identify $\PP$ with $\PP^\perp$ and $\QQ$ with $\QQ^\perp$
by the isomorphism $x\mapsto\omega(x,\cdot)$ between $\TT$ and $\EE$.
So, we will allow us to write $\EE=\PP\oplus\QQ$.

The connection $\nabla$ gives a derivation of the Fedosov algebra
$\Wb$, which is an extension
of the de Rham differential on functions.
Analogously, $\nabla$ gives such derivations of the algebras $T(\EE)\ot\La(\EE)$,
$S(\EE)\ot\La(\EE)$, and $(S(\PP)\ot S(\QQ))\ot\La(\EE)$.
These derivations commute with the maps
(\ref{rels}) and $\wi_{\PP,\QQ}$.

For convenience, we will mark the elements of the Fedosov algebra
lying in $\EE\ot 1$ by letters with hat over them ($\hat{x}$),
while by $dx$ we will
denote the copy of $\hat{x}$ lying in $1\ot\La^1\EE$.

Let $\omega=\omega_{ij}dx_i\wedge dx_j$ in some local coordinates.
It is easy to check that for
$\tilde{\delta}=\omega_{ij}\hat{x}_i\ot dx_j$ one has
\begin{equation}\label{reldelta}
\begin{split}
\delta&=\frac{1}{t}\ad(\tilde{\delta}), \\
\tilde{\delta}^2&=t\omega.
\end{split}
\end{equation}
Since the torsion of $\nabla$ is equal to zero,
\be{}\label{torsion}
\nabla(\tilde\delta)=0.
\ee{}
Since $\nabla^2$ is a $\CM$-linear derivation of degree $0$ preserving $\PP$,
there is an element $R\in \rho_{\PP,\EE}(\frak{sp}(\PP,\EE))\ot\La^2(\EE)$
such that
$\nabla^2=\frac{1}{t}\ad(R)$. In particular,
one has to be
\be{}\label{FPR}
F^\PP(R)\geq 1.
\ee{}
According to Fedosov, \cite{Fe}, we also define
$R^F\in\rho_\EE({\frak{sp}(\PP,\EE)})\ot\La^2(\EE)$
satisfying $\nabla^2=\ft\ad(R^F)$.
It follows from (\ref{torsion})
\be{}
\delta(R)=\delta(R^F)=0.
\ee{}
Following to Fedosov, we will consider connections on $\Wb$ of the form
\be{}
D=\nabla+\ft\ad(\gamma), \quad\quad \gamma\in \WW\ot\La^1(\EE).
\ee{}
We define the Wick curvature of $D$ as
\be{*}
\Omega_D=R+\nabla(\gamma)+\ft\gamma^2.
\ee{*}
According to Fedosov, we also define the Weyl (or Fedosov) curvature of $D$ as
\be{}\label{weylcur}
\Omega_D^F=R^F+\nabla(\gamma)+\ft\gamma^2.
\ee{}
Since by Lemma \ref{lem1.1}
$R^F=R+\frac{t}{2}tr(\nabla^2|_\PP)$,
we have
\be{}\label{wickcur}
\Omega^F_D=\Omega_D+\frac{t}{2}tr(\nabla^2|_\PP).
\ee{}
One checks
\be{}\label{flatness}
D^2=\ft\ad(\Omega_D)=\ft\ad(\Omega_D^F).
\ee{}
Let us take $\gamma$ in the form
\be{}
\gamma=\tilde\delta+r, \quad\quad r\in \WW\ot\La^1(\EE), \quad F^T(r)\geq 3.
\ee{}
Then the connection $D$ has the form
\be{}\label{Fcon}
\nabla+\delta+\ft\ad(r).
\ee{}
Using (\ref{reldelta}) and (\ref{torsion}), we obtain that its
Wick curvature is
\be{}\label{Fcurv}
\Omega_D=R+\nabla(\tilde\delta+r)+\frac{1}{t}(\tilde\delta+r)^2=
\omega+\delta r+R+\nabla r+\ft r^2.
\ee{}

\begin{propn}\label{prop1.3}
There exists an element $r\in \WW(\EE)\ot\La^1(\EE)$ such that

a) $F^T(r)\geq 3$;

b) $F^\PP(r)\geq 1$;

c) the connection $D=\nabla+\delta+\ft\ad(r)$ is flat, i.e.
$D^2=0$;

d) for its Wick curvature one has
$$\Omega_D=\omega;$$

e) its Weyl curvature $\Omega_D^F$ belong to $\Gamma(M,d\PP^\perp)$
and there is the formula
\be{}\label{Wc}
\Omega_D^F=\omega+\ftt tr(\nabla^2|_\PP).
\ee{}
\end{propn}

\begin{proof}
First of all, we apply the Fedosov method, (\cite{Fe}, Theorem 5.2.2),
to find $r$ satisfying d).
According to (\ref{Fcurv}), $r$ must obey  the
equation
\be{}
\delta r=-(R+\nabla r+\ft r^2)
\ee{}
Let us look for $r$ being the limit of the sequence, $r=\lim r_k$,
where $r_k\in \WW(\EE)\ot\La^1(\EE)$, $k=3,4,...$,
and $F^T(r_{k}-r_{k-1})\geq k$.
As in Lemma 5.2.3 of \cite{Fe},
using Proposition \ref{prop1.2} d) and the fact that $F^T(R)\geq 2$,
such $r_k$ can be calculated recursively:
\begin{equation}\label{calr}
\begin{split}
r_3&=-\delta_{\PP,\QQ}^{-1}(R) \\
r_{k+3}&=-(r_3+\delta_{\PP,\QQ}^{-1}(\nabla  r_{k+2}+
\frac{1}{t}r_{k+2}^2)).
\end{split}
\end{equation}
So, a) and d) are proven.

Let us prove that $F^\PP(r_k)\geq 1$ for all $k$.
The inequality $F^\PP(r_3)\geq 1$  follows from the fact
that $F^\PP(R)\geq 1$ and from Proposition (\ref{prop1.2}) c).
Suppose that $F^\PP(r_i)\geq 1$ for $i<k+3$, $k>0$. Then
$F^\PP(\nabla r_{k+2})\geq 1$ because $\nabla$ preserves $\PP$.
On the other hand, $F^\PP(r_{k+2}^2)\geq 1$
because of Proposition (\ref{prop1.2}) b),
therefore, as follows from (\ref{calr}),
$F^\PP(r_{k+3})\geq 1$ as well.
So, we have that $r$ being the limit of the convergent sequence
$r_k$ satisfies the conditions a), b), and d) of
the proposition.

c) obviously follows from d) and  (\ref{flatness}).

e) follows immediately from d), (\ref{wickcur}), and
Lemma \ref{lempolform}.
\end{proof}

\begin{propn}\label{prop1.3a}
Let $\tilde\nabla$ be another $\PP$-symplectic connection on $M$.
Let $\tilde r\in\WW(\EE)\ot\La^1(\EE)$, satisfy the conclusions of
Proposition \ref{prop1.3}, in particular,
the connection $\tilde D=\tilde\nabla+\delta+\ft\ad(\tilde r)$
is flat.
Then, there exists an element $B\in\WW(\EE)$ such that

a) $F^T(B)\geq 3$;

b) $F^\PP(B)\geq 1$

c) $e^{\ft\ad B}D=\tilde D$.
\end{propn}

\begin{proof}
Note that $\nabla-\tilde\nabla$ can be presented as $\ft\ad R_0$,
where $R_0\in\rho_{\PP,\EE}(\frak{sp}(\PP,\EE))$. Therefore,
$F^T(R_0)\geq 2$ and $F^\PP(R_0)\geq 1$.
Let us put $R_1=r-\tilde r$. Then,
$F^T(R_1)\geq 3$ and $F^\PP(R_1)\geq 1$.
We have
$$\tilde D=D-\ft\ad(R_0+R_1).$$
Since $\Omega_D=\Omega_{\tilde D}=\omega$,
using (\ref{torsion}) we obtain
$$\delta(R_0)=0.$$
It follows that the element $B_0=\delta_{\PP,\QQ}^{-1}(R_0)$
is such that $\delta(B_0)=R_0$ and
$F^T(B_0)\geq 3$, $F^\PP(B_0)\geq 1$.

Replacing $D$ by $e^{\ft\ad B_0}D$ we obtain
$$\tilde D=D-\ft\ad(R'_0+R'_1),$$
where $F^T(R'_0)\geq 3$, $F^\PP(R'_0)\geq 1$ and
$F^T(R'_1)\geq 4$, $F^\PP(R'_1)\geq 1$.

Proceeding by induction on $F^T$-filtration, we obtain
a sequence $B_i\in\WW(\EE)$ with increasing $F^T$-filtration
and such that $\Pi_{i=0}^\infty e^{\ft\ad B_i}(D)=\tilde D$.

Since elements $e^{\ft\ad B'}$, $B'\in\WW(\EE)$, $F^T(B)\geq 3$,
form a pro-unipotent Lie group,
there exists
an element $B\in\WW(\EE)$, $F^T(B)\geq 3$, $F^\PP(B)\geq 1$,
such that $\Pi_{i=0}^\infty e^{\ft\ad B_i}=e^{\ft\ad B}$.
\end{proof}

Let $D$ be a connection satisfying Proposition \ref{prop1.3}.
Denote by $\WW_D$ the subsheaf of $\WW$ consisting of
flat sections $a$, i.e. such that $Da=0$.
Since $D$ is a derivation of $\Wb$, it is clear that $\WW_D$ is a
sheaf of subalgebras.
Let $\sigma=id-(\delta\delta_{\PP,\QQ}^{-1}+\delta_{\PP,\QQ}^{-1}\delta)$.
Then, as follows from (\ref{relde}), $\sigma: \WW\to\CM$ is a projection,
where $\CM$ is considered as the center of the algebra $\WW$.

\begin{propn}\label{prop1.4}
a) The map $\sigma:\WW_D\to\CM$ is a bijection.

b) The inverse map $\eta:\CM\to \WW_D$ has the form
$\eta(f)=f+\hat f$, there $F^T(\hat f)>F^T(f)$.

c) If $df\in\PP$, then $F^\PP(\hat f)\geq 1$.

d) If $df\in\PP$, then $\sigma(\eta(f)\eta(g))=fg$
for any $g\in\CM$.
\end{propn}

\begin{proof}
Again, we apply the Fedosov iteration procedure.
According to \cite{Fe}, Theorem 5.2.4, we look for $\eta(f)$
as a limit, $\eta(f)=\lim a_k$, there $a_k\in \WW$ can
be calculated recursively:
\begin{equation}\label{calf}
\begin{split}
a_0&=f  \\
a_{k+1}&=a_0+\delta_{\PP,\QQ}^{-1}(\nabla a_k+\ft\ad r(a_k)).
\end{split}
\end{equation}
Put $\hat f=\eta(f)-a_0$.
As in \cite{Fe}, Theorem 5.2.4, one proves that such $\eta(f)$
and $\hat f$ satisfy a) and b).
Now observe that $a_1-a_0=\delta_{\PP,\QQ}^{-1}(1\ot df)$ and
if  $df\in\PP$, then $F^\PP(a_1-a_0)\geq 1$.
By induction, we conclude that  $F^\PP(a_k-a_0)\geq 1$ for all $k\geq 1$.
So $F^\PP(a-a_0)\geq 1$ as well, which proves c).

Let us prove d). We have $\eta(f)\eta(g)=f\eta(g)+\hat f\eta(g)$.
Since by c) $F^\PP(\hat f)\geq 1$, $F^\PP(\hat f\eta(g))\geq 1$ as well.
It follows that $\sigma(\hat f\eta(g))=0$ and
$\sigma(\eta(f)\eta(g))=\sigma(f\eta(g))=fg$,
because $\sigma$ is a $\CM$-linear map and $\sigma(\eta(g))=g$.
\end{proof}

\subsection{Existence of PSP's}

Let $(M,\omega,\PP)$ be a polarized symplectic manifold.
Recall that in Subsection \ref{subsclass} we have assigned
to any PSP $(\mu_t, \OO_t)$ on $(M,\omega,\PP)$ an element
$\tau(\mu_t,\OO_t)\in\YY(M,\omega,\PP)$, which is a pair
$(\omega_t,\PP_t)$
being a deformation of the pair $(\omega,\PP)$.
The form $\omega_t$ represents the extension class
$\cle(\mu_t,\OO_t)$.

We show now that any element of $\YY(M,\omega,\PP)$
corresponds to a PSP.

\begin{propn}\label{essent2}
a) For any pair $(\omega_t,\PP_t)\in\YY(M,\omega,\PP)$,
there exists a PSP, $(\mu_t,\OO_t)$, such that
\be{}\label{lpa}
\tau(\mu_t,\OO_t)=(\omega_t,\PP_t).
\ee{}
b) The Fedosov class of the corresponding star-product
$(\CM,\mu_t)$ is represented by the form of
$\omega+t\Gamma(M,d(\PP_t^\perp))$ equal to
\be{}\label{lpb}
\theta_t=\omega_t+\ftt tr(\nabla^2|_{\PP_t}),
\ee{}
where $\nabla$ is a $\PP_t$-symplectic connection on
the formal symplectic manifold $(M,\omega_t,\PP_t)$.
\end{propn}

\begin{proof}
Let $\Wb$ be the Fedosov
algebra on $M$ corresponding to the symplectic form $\omega_t$.
Let $\nabla$ be a $\PP_t$-symplectic connection on $M$ corresponding to
$\omega_t$ and $D$ the flat connection on $\Wb$ constructed in
Proposition \ref{prop1.3} c). Let $\WW_D$ be the sheaf of
flat sections of $\WW$. Define a star-product $(\CM,\mu_t)$ on $M$
carrying over
the multiplication from $\WW_D$ to $\CM$ via the map $\sigma$ from
Proposition \ref{prop1.4}. Point d) of that proposition shows
that, in fact, this star-product present
the PSP $(\CM,\mu_t,\OO_{\PP_t})$.
We are going to prove that this star-product is as required.
In the following we identify $\WW_D$ with the
corresponding PSP via $\sigma$.

Let us prove (\ref{lpa}).
Let $(U_\alpha)$ be an open covering of $M$
such that on each $U_\alpha$ there exist
formal Darboux
coordinates $x_{\alpha i}$, $y_{\alpha i}$,
$x_{\alpha i}\in\OO_{\PP_t}$ with respect to $\omega_t$.

Denote by $\nabla_\alpha$ the standard flat $\PP_t$-symplectic
connections over $U_\alpha$ such that the forms $dx_{\alpha i}$,
$dy_{\alpha i}$ are flat sections in $\PP^\perp_t$.
Then, the connections $D_\alpha=\nabla_\alpha+\delta$
satisfy on $U_\alpha$ Proposition \ref{prop1.3} with $r=0$.
Let $\WW_{D_\alpha}$ be the star-product on $U_\alpha$
constructed in Proposition \ref{prop1.4} via
flat sections of $D_\alpha$. It is easy to see that
$\WW_{D_\alpha}$ coincides with the Moyal-Wick
PSP with respect
to the Darboux coordinates $x_{\alpha i}$, $y_{\alpha i}$
(see Example \ref{examp2}),
so $x_{\alpha i}$, $y_{\alpha i}$ are also Darboux coordinates
for the bracket $\ft[\cdot,\cdot]$ in $\WW_{D_\alpha}$.

Since $D$ and $D_\alpha$ have the same Wick curvature $\omega_t$,
there exist, by Proposition \ref{prop1.3a}, elements $B_\alpha\in\WW$
such that $e^{\ft\ad B_\alpha}D_\alpha=D$.
It is clear that $e^{\ft\ad B_\alpha}$ acting on $\WW$
takes $\WW_{D_\alpha}$ to $\WW_D$,
and point b) of that
Proposition implies that it is identical on $\OO_{\PP_t}$.

Let $\Psi_{\alpha,\beta}=e^{-\ft\ad B_\alpha}e^{\ft\ad B_\beta}$.
These may be considered as isomorphisms over $U_\alpha\cap U_\beta$
gluing star-products $\WW_{D_\alpha}$ to a global star-product on $M$
that is obviously isomorphic to $\WW_D$.
We see, that functions $x_{\alpha i}$, $y_{\alpha i}$ form
local Darboux coordinates corresponding to that star-product.
So, the characteristic 2-form $\clP(\WW_D)$
(see Subsection \ref{subsclass}) is locally represented as
$dy_{\alpha i}\wedge dx_{\alpha i}$.
On the other hand, this form is equal to $\omega_t$, since
from very beginning the
functions $x_{\alpha i}$, $y_{\alpha i}$ have been chosen as
Darboux coordinates for it. Hence, $\tau(\WW_D)=(\omega_t,\PP_t)$.

b) Follows from Proposition \ref{prop1.3} e) and Lemma \ref{lempolform}.
\end{proof}

\section{The main theorem and corollaries}

Let $(M,\omega,\PP)$ be a polarized symplectic manifold.
Denote by $\YY$ the set of pairs $(\omega_t,\PP_t)$, where
$\omega_t=\omega+t\omega_1+\cdots$ is a deformed symplectic form
and $\PP_t$, $\PP_0=\PP$, its polarization. Let $Aut(M)$ be
the group of formal automorphisms of $M$.

\begin{thm}\label{thmm}
a) The equivalence classes of PDQ's on $(M,\omega,\PP)$
are in one-to-one correspondence with the orbits in $\YY$ under
the $Aut(M)$-action.

b) Let the pair $(\omega_t,\PP_t)$ be a point on the orbit
corresponding to a PDQ $(\bA_t,\bO_t)$.
Then, $(\bA_t,\bO_t)$ is isomorphic to a
PSP, $(\CM,\mu_t,\OO_t)$, where $\OO_t$
consists of functions constant along $\PP_t$ and the multiplication
$\mu_t$ satisfies the condition
\be{*}
\mu_t(f,g)=fg \quad\quad \mbox{for}\ \ f\in\OO_t,\ g\in\CM.
\ee{*}

c) The form $\omega_t$ represents the extension class
$\cle(\bA_t,\bO_t)\in H^1(M,\Omega^{1cl}_{\OO_t})$,
associated with $(\bA_t,\bO_t)$.

d) Under the hypothesis of b),
the Fedosov class of the deformation quantization $\bA_t$
can be represented by a form $\theta_t$ that is a deformation
of $\omega$ and polarized by $\PP_t$.
It is defined by the formula
\be{}\label{eqFC}
\theta_t=\omega_t+\ftt tr(\nabla^2|_{\PP_t}),
\ee{}
where $\nabla$ is a $\PP_t$-symplectic connection on
the formal symplectic manifold $(M,\omega_t,\PP_t)$.
\end{thm}

\begin{proof}
Parts a) and b) follow from Proposition \ref{eqcl},
Corollary \ref{corin}, and
Proposition \ref{essent2} a).
Part c) follows from Lemma \ref{lemrepcl}.
Part d) is the same as Proposition \ref{essent2} b).
\end{proof}

\begin{remark}\label{remPC}
We have interpreted the form $\omega_t$ from (\ref{eqFC}) as a
representative of the extension class associated to a PDQ (see
Section 8). The form $tr(\nabla^2|_{\PP_t})$ is the curvature form
of the connection induced by $\nabla|_{\PP_t}$ on the complex line
bundle $\det(\PP_t)$. Actually, $\det(\PP_t)$ can be presented as
a line $\OO_t$-bundle, i.e. as a locally free sheaf of $\OO_t$-modules
of rank one. Indeed, $\det(\PP_t)=\CM\ot_{\OO_t}\LL$, where
$\LL=\Omega_{\OO_t}^n$, $n=\half\dim M$.
The form
$-tr(\nabla^2|_{\PP_t})$, as well as {\em minus} curvature form of
any other connection on $\LL$, can be interpreted as a representative of
the extension class of an $\OO_t$-extension of $Der(\OO_t)$ associated
with $\LL$. Indeed, let $\widetilde T_\LL$ denote the sheaf of
$\OO_t$-differential operators on $\LL$ of order at most one. Then
$\widetilde T_\LL$ equipped with the left $\OO_t$-module structure
and the Lie bracket given by the commutator naturally forms a
$\OO_t$-extension of $Der(\OO_t)$. Splittings of this extension are
flat connection on $\LL$. Let $d_\alpha$ be local flat connections
on $\LL$ in some open covering $\{U_\alpha\}$ of $M$. Then,
$d_\alpha-d_\beta$ are closed 1-forms of $\Omega_{\OO_t}^1$
that form a \v Cech cocycle.
Hence, there exist smooth 1-forms $f_\alpha\in\PP^\perp_t$
such that
$d_\alpha-d_\beta=f_\alpha-f_\beta$. Differential operators
$d_\alpha-f_\alpha$ form a global connection on $\LL$,
$\nabla_\LL$, with the curvature locally equal to $-df_\alpha$. On
the other hand, by definition (see Section 8),
the extension class of $\widetilde T_\LL$ is
represented by the form $df_\alpha\in d\PP^\perp_t$,
i.e. $-\nabla_\LL^2$ .
So, projecting the equality (\ref{eqFC}) to
$\Hwt$, we obtain
\be{}
[\theta_t]=\cle(\mu_t,\OO_t)-\ftt\cle(\widetilde T_{\det(\PP_t)}).
\ee{}
Details are left to the reader.

Note that the form $-\frac{1}{2\pi\sqrt{-1}}tr(\nabla^2|_{\PP_t})$
represents the first Chern class of $\PP$, \cite{KN}.
\end{remark}

\begin{cor}
Let $\bA_t$ be a
deformation quantization on $(M,\omega)$. Suppose its
Fedosov class $cl_F(\bA_t)$ is represented by the
form $\theta_t$ that has a polarization $\PP_t$.
Then $\bA_t$ can be extended to a PDQ $(\bA_t,\bO_t)$, where
$\bO_t$ is isomorphic to $\OO_{\PP_t}$.
\end{cor}
\begin{proof}
Let $\nabla$ be a $\PP_t$-symplectic connection
on the formal symplectic manifold $(M,\theta_t,\PP_t)$.
Let $(\CM,\mu_t,\OO_t)$ be a PSP such that
$$\tau(\CM,\mu_t,\OO_t)=(\theta_t-\ftt tr(\nabla^2|_{\PP_t}),\PP_t).$$
By (\ref{eqFC}), $cl_F(\CM,\mu_t)=[\theta_t]=cl_F(\bA_t)$,
therefore star-products $\bA_t$ and $(\CM,\mu_t)$ are
equivalent.
\end{proof}

\begin{remark} All constructions of the paper can be extended to
the case when $M$ is a formal manifold,
$M_\lambda$, which is $M$ endowed with the function sheaf
$C_M^\infty[[\lambda]]$, $\lambda$ a formal parameter.
A formal polarized symplectic manifold is a triple,
$(M_\lambda,\omega_\lambda,\PP_\lambda)$, where $\omega_\lambda$
is a formal symplectic form on $M_\lambda$ and $\PP_\lambda$ its polarization.
The above construction of a polarized star-product applied to a formal
polarized symplectic manifold $(M_\lambda,\omega_\lambda,\PP_\lambda)$
gives the following
\begin{propn}
Let $(\bA_t,\bO_t)$ be a PDQ and $(\omega_t,\PP_t)$ an element on
the orbit corresponding to $(\bA_t,\bO_t)$.
Then, there exists on $(M_\lambda,\omega_\lambda,\PP_\lambda)$
a PSP
$$(C_M^\infty[[\lambda]][[t]],\mu_{\lambda,t},\OO_{\lambda,t})$$
such that $(\bA_t,\bO_t)$ is equivalent to the diagonal sub-family
$(C_M^\infty[[t]],\mu_{t,t},\OO_{t,t})$ obtained by the substitution
$\lambda=t$.
\end{propn}
\end{remark}

\bigskip

e-mail: donin@macs.biu.ac.il
\end{document}